# ADAPTIVE GOODNESS-OF-FIT TESTS IN A DENSITY MODEL

BY MAGALIE FROMONT AND BÉATRICE LAURENT

*Université Rennes II and INSA de Toulouse*

Given an i.i.d. sample drawn from a density $f$, we propose to test that $f$ equals some prescribed density $f_0$ or that $f$ belongs to some translation/scale family. We introduce a multiple testing procedure based on an estimation of the $\mathbb{L}_2$-distance between $f$ and $f_0$ or between $f$ and the parametric family that we consider. For each sample size $n$, our test has level of significance $\alpha$. In the case of simple hypotheses, we prove that our test is adaptive: it achieves the optimal rates of testing established by Ingster [*J. Math. Sci.* **99** (2000) 1110–1119] over various classes of smooth functions simultaneously. As for composite hypotheses, we obtain similar results up to a logarithmic factor. We carry out a simulation study to compare our procedures with the Kolmogorov–Smirnov tests, or with goodness-of-fit tests proposed by Bickel and Ritov [in *Nonparametric Statistics and Related Topics* (1992) 51–57] and by Kallenberg and Ledwina [*Ann. Statist.* **23** (1995) 1594–1608].

**1. Introduction.** Suppose that we observe $n$ independent and identically distributed (i.i.d.) real random variables $X_1, \ldots, X_n$ with common unknown density $f$. Let $f_0$ be some specified density. In this paper we consider the problem of testing the null hypothesis "$f \in \mathcal{F}$" against "$f \notin \mathcal{F}$" where $\mathcal{F}$ equals either the singleton $\{f_0\}$ or the parametric family

$$(1.1) \qquad \mathcal{F} = \left\{ \frac{1}{\sigma} f_0 \left( \frac{\cdot - \mu}{\sigma} \right), (\mu, \sigma) \in K \right\},$$

for some subset $K$ of $\mathbb{R} \times ]0, +\infty[$.

This problem has been widely studied since the famous Kolmogorov–Smirnov and Cramér–von Mises tests based on the empirical distribution function.

Assuming that $f$ belongs to $\mathbb{L}_2(\mathbb{R})$, it is quite natural to construct a test based on the estimation of the squared $\mathbb{L}_2$-distance between $f$ and $\mathcal{F}$. In order to test the simple hypothesis "$f = f_0$," we actually consider a suitable collection of estimators of $\int_{\mathbb{R}}(f - f_0)^2$ and decide to reject the null hypothesis if some estimator in the collection is larger than its $(1 - u_\alpha)$ quantile under the null hypothesis, $u_\alpha$ being calibrated so that the final test has level of significance $\alpha$. We then generalize this procedure to test that $f$ belongs to the translation/scale family given by (1.1).

From a theoretical point of view, we evaluate the performances of our tests in terms of uniform separation rates with respect to the $\mathbb{L}_2$-distance over classes of smooth functions. Given $\beta$ in $]0, 1[$ and a class of functions $\mathcal{B} \subset \mathbb{L}_2(\mathbb{R})$, we define the uniform separation rate $\rho(\Phi_\alpha, \mathcal{B}, \beta)$ of a level $\alpha$ test $\Phi_\alpha$ of the null hypothesis "$f \in \mathcal{F}$" over the class $\mathcal{B}$ as the smallest number $\rho$ such that the test guarantees a power at least equal to $(1 - \beta)$ for all alternatives $f$ in $\mathcal{B}$ at a distance $\rho$ from $\mathcal{F}$. More precisely, denoting by $d_2(f, \mathcal{F})$ the $\mathbb{L}_2$-distance between $f$ and $\mathcal{F}$ and by $\mathbb{P}_f$ the distribution of the observation $(X_1, \ldots, X_n)$,

$$\rho(\Phi_\alpha, \mathcal{B}, \beta)$$
$$= \inf\{\rho > 0, \ \forall f \in \mathcal{B}, \ d_2(f, \mathcal{F}) \geq \rho \Rightarrow \mathbb{P}_f(\Phi_\alpha \text{ rejects}) \geq 1 - \beta\}.$$

Assuming that $f$ belongs to $\mathcal{B}$, the uniform separation rate $\rho(\Phi_\alpha, \mathcal{B}, \beta)$ is asymptotically related to the minimax rate of testing $\rho_n$ introduced by Ingster [14] and referred to as the *critical radius*. Indeed, by definition, $\rho_n \to 0$ as $n \to +\infty$ and satisfies:

(a) For any sequence $\rho'_n$ such that $\rho'_n/\rho_n = o_n(1)$,

$$\inf_{\Phi_n}\left\{\sup_{f \in \mathcal{F}} \mathbb{P}_f(\Phi_n = 1) + \sup_{f \in \mathcal{B}, d_2(f, \mathcal{F}) \geq \rho'_n} \mathbb{P}_f(\Phi_n = 0)\right\} = 1 - o_n(1),$$

where the infimum is taken over all tests $\Phi_n$ with values in $\{0, 1\}$ rejecting the null hypothesis "$f \in \mathcal{F}$" when $\Phi_n = 1$.

(b) For any $\alpha$, $\beta > 0$, there exist some constant $C > 0$ and some test $\Phi_n^*$ such that the two following inequalities hold:

$$(1.2) \qquad\qquad \sup_{f \in \mathcal{F}} \mathbb{P}_f(\Phi_n^* = 1) \leq \alpha + o_n(1),$$

$$(1.3) \qquad \sup_{f \in \mathcal{B}, d_2(f, \mathcal{F}) \geq C\rho_n} \mathbb{P}_f(\Phi_n^* = 0) \leq \beta + o_n(1).$$

Since the goodness-of-fit test to some specified density $f_0$ can be reduced to a test of uniformity on $[0, 1]$ for the variables $F_0(X_i)$ (where $F_0$ is the distribution function associated with the density $f_0$), many papers are devoted to the problem of testing uniformity on $[0, 1]$. The main reference for the computation of minimax rates of testing for this problem is the series of papers due to Ingster [14]. In particular, under the prior assumption that $f$ belongs to some Hölder class with smoothness parameter $s > 0$, Ingster establishes the minimax rate of testing $\rho_n = n^{-2s/(4s+1)}$. But the tests proposed to ensure the achievement of this rate [namely the inequalities given in (1.2) and (1.3)] are structurally based on the prior assumption; this is a crucial problem for their practical application since the smoothness parameter $s$ is typically unknown. Following the work of Spokoiny [23] in the Gaussian white noise model, Ingster [15] focuses on the problem of finding an adaptive (assumption-free) test of uniformity on $[0,1]$. He proves that adaptation is not possible without some loss of efficiency of the order of an extra $\log \log n$ factor and he presents an adaptive test which is based on chi-square statistics.

Other methods having Neyman's test as starting point are proposed in order to avoid using any prior assumption on the smoothness of $f$. To test uniformity on $[0,1]$, Neyman [20] suggests considering some orthonormal basis $\{\phi_l, \ l \geq 0\}$ of $\mathbb{L}_2([0,1])$ with $\phi_0 = \mathbb{I}_{[0,1]}$ and rejecting the null hypothesis "$f = \mathbb{I}_{[0,1]}$" if the estimator $\sum_{l=1}^{D}(\sum_{i=1}^{n}\phi_l(X_i)/n)^2$ of $\theta_D = \sum_{l=1}^{D}(\int_{[0,1]} f\phi_l)^2$ is large enough, where $D$ is some given integer. Bickel and Ritov [4], Ledwina [19] and Kallenberg and Ledwina [17] introduce data-driven versions of Neyman's test where the parameter $D$ is chosen via some penalized criterion. Inglot and Ledwina [13] establish theoretical results for the test described in Kallenberg and Ledwina [17]. These results which essentially deal with the asymptotic efficiency of the test with respect to the Neyman–Pearson test do not, however, lead to any optimality of the uniform separation rates. Fan [8] also proposes a new version of Neyman's test based on wavelet thresholding to test that the mean of a Gaussian vector equals 0 with applications to goodness-of-fit tests in a density model. When we test uniformity on $[0,1]$, our method amounts to considering, for all integer $D$ in some set $\mathcal{D}_n$, the unbiased estimator of $\theta_D$ defined by

$$\hat{\theta}_D = \frac{1}{n(n-1)} \sum_{l=1}^{D} \sum_{\substack{i \neq j=1}}^{n} \phi_l(X_i)\phi_l(X_j)$$

and to penalizing this estimator by its $(1 - u_\alpha)$ quantile under the null hypothesis. The main difference between our method and the testing procedures proposed by previous authors lies in the order of magnitude of the penalty term. While Ledwina [19] and Kallenberg and Ledwina [17] choose the parameter $D$ by using Schwarz's Bayesian information criterion (BIC), Kallenberg [16] gives a discussion of the choice of the penalties for data-driven Neyman's tests. But the criteria considered in these papers have been introduced to estimate the density $f$ itself, whereas our penalties correspond to the ones used to build adaptive estimators of $\int_{\mathbb{R}} f^2$ by model selection in [18]. This choice allows us to obtain optimal uniform separation rates with respect to the $\mathbb{L}_2$-distance.

As for testing a composite null hypothesis, Pouet [22] proves that provided that $f$ belongs to $\mathbb{L}_2([0,1])$ and some Hölder class, the minimax rate of testing is comparable to the rate for the simple hypothesis "$f = \mathbb{I}_{[0,1]}$." However, the test proposed by Pouet depends on the smoothness assumption on $f$, which is not satisfactory from an experimental point of view. Inglot, Kallenberg and Ledwina [12] introduce a procedure using no prior information about the smoothness of $f$ to test composite hypotheses like "$f \in \{f(x,\beta), \ \beta \in B\}$" with $B \subset \mathbb{R}^q$. This procedure, generalizing Kallenberg and Ledwina's one [17], also consists of a combination of Neyman's smooth test and Schwarz's selection rule. Its construction is based on the consideration of a sequence of exponential families with increasing dimensions to describe departures from the null model. The "right" dimension is selected by an extended Schwarz's rule, which is obtained by inserting the maximum likelihood estimator $\hat{\beta}$ of $\beta$ under the null hypothesis in the original definition of Schwarz's BIC. The next step is the application of Neyman's smooth test using a quadratic score statistic in the selected dimension. Inglot, Kallenberg and Ledwina [12] prove the consistency of the test at essentially any alternative.

The approach considered in the present paper has been initiated by Baraud, Huet and Laurent [1, 2, 3] for the problem of testing linear or qualitative hypotheses in the Gaussian regression model. The properties of the testing procedures proposed here are nonasymptotic. For each $n$, the tests have the desired level of significance and we characterize some sets of alternatives over which they have a prescribed power. For the problem of testing goodness-of-fit of some specified density, we state in Section 2 that our procedure is adaptive over some collection of classes of smooth functions in the sense that it achieves the optimal "adaptive" rate of testing established by Ingster [15] over all the classes of the collection simultaneously. We also investigate in Section 4 the test from a practical point of view by Monte Carlo experiments. The results show that our procedure is competitive with the ones due to Bickel and Ritov [4] or Kallenberg and Ledwina [17]. For the problem of testing the hypothesis "$f \in \mathcal{F}$," where $\mathcal{F}$ is the translation/scale family defined by (1.1), we get in Section 3 uniform separation rates over classes of smooth alternatives of the same order (up to a logarithmic factor) as the rates obtained when testing the simple hypothesis "$f = f_0$." We finally implement the procedure to test exponentiality in Section 4; we can notice that it gives particularly good results in comparison with the Kolmogorov–Smirnov test under oscillating alternatives. The proofs of the results stated in the paper are detailed in Section 5.

**2. A goodness-of-fit test.** Let $X_1, \ldots, X_n$ be i.i.d. random variables with common density $f$ with respect to the Lebesgue measure on $\mathbb{R}$. Let $f_0$ be some given density in $\mathbb{L}_2(\mathbb{R})$ and let $\alpha$ be in $]0, 1[$. Assuming that $f$ belongs to $\mathbb{L}_2(\mathbb{R})$, we construct a level $\alpha$ test of the null hypothesis "$f = f_0$" against the alternative "$f \neq f_0$" from the observation $(X_1, \ldots, X_n)$.

In the following, $\|\cdot\|_2$ and $\langle\cdot, \cdot\rangle$, respectively, denote the usual norm and scalar product in $\mathbb{L}_2(\mathbb{R})$. For any bounded function $g$, $\|g\|_\infty = \sup_{x \in \mathbb{R}} |g(x)|$.

2.1. *Description of the test.* Our test is based on an estimation of the quantity $\|f - f_0\|_2^2$ that is $\|f\|_2^2 + \|f_0\|_2^2 - 2\langle f, f_0 \rangle$. Since $\langle f, f_0 \rangle$ is usually estimated by the empirical estimator $\sum_{i=1}^n f_0(X_i)/n$, the key point is the estimation of $\|f\|_2^2$. As in [18], we introduce an at most countable collection $\{S_m, \ m \in \mathcal{M}\}$ of linear subspaces of $\mathbb{L}_2(\mathbb{R})$. For all $m$ in $\mathcal{M}$, let $\{p_l, l \in \mathcal{L}_m\}$ be some orthonormal basis of $S_m$. The variable

$$(2.1) \qquad \hat{\theta}_m = \frac{1}{n(n-1)} \sum_{l \in \mathcal{L}_m} \sum_{i \neq j = 1}^n p_l(X_i) p_l(X_j)$$

is an unbiased estimator of $\|\Pi_{S_m}(f)\|_2^2$, where $\Pi_{S_m}$ denotes the orthogonal projection onto $S_m$. Then $\|f - f_0\|_2^2$ can be estimated by

$$\hat{T}_m = \hat{\theta}_m + \|f_0\|_2^2 - \frac{2}{n} \sum_{i=1}^n f_0(X_i),$$

for any $m$ in $\mathcal{M}$. Denoting by $t_m(u)$ the $(1 - u)$ quantile of the law of $\hat{T}_m$ under the hypothesis "$f = f_0$" and considering

$$(2.2) \qquad u_\alpha = \sup\left\{ u \in ]0, 1[, \ \mathbb{P}_{f_0}\left( \sup_{m \in \mathcal{M}} (\hat{T}_m - t_m(u)) > 0 \right) \leq \alpha \right\},$$

we introduce the test statistic $T_\alpha$ defined by

$$T_\alpha = \sup_{m \in \mathcal{M}} (\hat{T}_m - t_m(u_\alpha)).$$

Our test consists of rejecting the null hypothesis if $T_\alpha$ is positive.

In practice, the values of $u_\alpha$ and the quantiles $\{t_m(u_\alpha), m \in \mathcal{M}\}$ are estimated by Monte Carlo experiments under $f_0$ as explained in Section 4.

This method amounts to a multiple testing procedure. Indeed, for each $m$ in $\mathcal{M}$, we construct a level $u_\alpha$ test of the null hypothesis "$f = f_0$" by rejecting this hypothesis if $\hat{T}_m$ is larger than its $(1 - u_\alpha)$ quantile under the hypothesis "$f = f_0$." We thus obtain a collection of tests and we decide to reject the null hypothesis if for some of the tests of the collection this hypothesis is rejected.

2.2. *The power of the test.* Let us now describe the collection of linear subspaces $\{S_m, m \in \mathcal{M}\}$ that we use to define our testing procedure here. This collection is obtained by mixing spaces generated by constant piecewise functions, scaling functions and, in the case of compactly supported densities, trigonometric polynomials.

(i) For all $D$ in $\mathbb{N}^*$ and $k$ in $\mathbb{Z}$, let

$$I_{D,k} = \sqrt{D} \mathbb{I}_{[k/D,(k+1)/D[}.$$

For all $D$ in $\mathbb{N}^*$, we define $S_{(1,D)}$ as the space generated by the functions $\{I_{D,k}, \ k \in \mathbb{Z}\}$ and

$$\hat{\theta}_{(1,D)} = \frac{1}{n(n-1)} \sum_{k \in \mathbb{Z}} \sum_{\substack{i \neq j=1}}^{n} I_{D,k}(X_i) I_{D,k}(X_j).$$

(ii) Let us consider a pair of compactly supported orthonormal wavelets $(\varphi, \psi)$ such that for all $J$ in $\mathbb{N}$, $\{\varphi_{J,k} = 2^{J/2}\varphi(2^J \cdot - k), k \in \mathbb{Z}\} \cup \{\psi_{j,k} = 2^{j/2}\psi(2^j \cdot - k), j \in \mathbb{N}, j \geq J, k \in \mathbb{Z}\}$ is an orthonormal basis of $\mathbb{L}_2(\mathbb{R})$. For all $J$ in $\mathbb{N}$ and $D = 2^J$, we define $S_{(2,D)}$ as the space generated by the scaling functions $\{\varphi_{J,k}, k \in \mathbb{Z}\}$ and

$$\hat{\theta}_{(2,D)} = \frac{1}{n(n-1)} \sum_{k \in \mathbb{Z}} \sum_{\substack{i \neq j=1}}^{n} \varphi_{J,k}(X_i) \varphi_{J,k}(X_j).$$

(iii) Let us consider the Fourier basis of $\mathbb{L}_2([0,1])$ given by

$$g_0(x) = \mathbb{I}_{[0,1]}(x),$$
$$g_{2p-1}(x) = \sqrt{2}\cos(2\pi p x)\mathbb{I}_{[0,1]}(x) \qquad \text{for all } p \geq 1,$$
$$g_{2p}(x) = \sqrt{2}\sin(2\pi p x)\mathbb{I}_{[0,1]}(x) \qquad \text{for all } p \geq 1.$$

For all $D$ in $\mathbb{N}^*$, we define $S_{(3,D)}$ as the space generated by the functions $\{g_l, l = 0, \ldots, D\}$ and

$$\hat{\theta}_{(3,D)} = \frac{1}{n(n-1)} \sum_{l=0}^{D} \sum_{\substack{i \neq j=1}}^{n} g_l(X_i) g_l(X_j).$$

We want here to notice that the constants involved in the following may depend on the chosen scaling function $\varphi$, but we will not always specify it.

Introduce $\mathbb{D}_1 = \mathbb{D}_3 = \mathbb{N}^*$ and $\mathbb{D}_2 = \{2^J, \ J \in \mathbb{N}\}$. For $l$ in $\{1, 2, 3\}$, $D$ in $\mathbb{D}_l$, $\Pi_{S_{(l,D)}}$ denotes the orthogonal projection onto $S_{(l,D)}$ in $\mathbb{L}_2(\mathbb{R})$.

For all $l$ in $\{1, 2, 3\}$, we take $\mathcal{D}_l \subset \mathbb{D}_l$ with $\bigcup_{l \in \{1,2,3\}} \mathcal{D}_l \neq \varnothing$ and $\mathcal{D}_3 = \varnothing$ if the $X_i$'s are not included in $[0,1]$. Let

$$\mathcal{M} = \{(l, D), l \in \{1, 2, 3\}, D \in \mathcal{D}_l\}.$$

For all $m$ in $\mathcal{M}$, we set

$$\hat{T}_m = \hat{\theta}_m + \|f_0\|_2^2 - \frac{2}{n} \sum_{i=1}^{n} f_0(X_i).$$

The test statistic that we consider is

$$(2.3) \qquad T_\alpha = \sup_{m \in \mathcal{M}} (\hat{T}_m - t_m(u_\alpha)),$$

where $t_m(u_\alpha)$ is defined in Section 2.1.

The aim of the following theorem is to describe classes of alternatives over which the corresponding test has a prescribed power.

THEOREM 1. *Let $X_1, \ldots, X_n$ be i.i.d. real random variables with common density $f$ and let $f_0$ be some given density. Let $T_\alpha$ be the test statistic defined by (2.3). Assume that $f_0$ and $f$ belong to $\mathbb{L}_\infty(\mathbb{R})$ and fix some $\beta$ in $]0, 1[$. For any $\varepsilon$ in $]0, 2[$, there exist some positive constants $C_1(\beta)$ and $C_2(\beta, \varepsilon, \|f\|_\infty, \|f_0\|_\infty)$ such that, setting for all $m = (l, D)$ in $\mathcal{M}$,*

$$V_m(\beta) = \frac{C_1(\beta)}{n} \left( (\sqrt{\|f\|_\infty} + \|f\|_\infty) \sqrt{D} + \frac{D}{n} \right) + \frac{C_2(\beta, \varepsilon, \|f\|_\infty, \|f_0\|_\infty)}{n},$$

*if $f$ satisfies*

$$\|f - f_0\|_2^2 > (1 + \varepsilon) \inf_{m \in \mathcal{M}} \{ \|f - \Pi_{S_m}(f)\|_2^2 + t_m(u_\alpha) + V_m(\beta) \},$$

*then*

$$\mathbb{P}_f(T_\alpha \le 0) \le \beta.$$

COMMENTS. (i) Let us see what is the advantage of considering a multiple testing procedure. We deduce from Theorem 1 that if we fix some element $m$ in $\mathcal{M}$ and if we focus on the test that rejects the null hypothesis when $\hat{T}_m$ is larger than its $(1 - \alpha)$ quantile under the hypothesis "$f = f_0$" denoted by $t_m(\alpha)$, then the error probability of the second kind of the test is smaller than $\beta$ for all $f$ such that

$$\|f - f_0\|_2^2 > (1 + \varepsilon) \{ \|f - \Pi_{S_m}(f)\|_2^2 + t_m(\alpha) + V_m(\beta) \}.$$

For the multiple testing procedure, the right-hand side of the above inequality is replaced by its infimum over all $m$ in $\mathcal{M}$, at the price that $t_m(\alpha)$ is replaced by $t_m(u_\alpha)$. When we evaluate the uniform separation rates, we show that for the collections $\{ S_m, \ m \in \mathcal{M} \}$ that we have chosen, the quantities $t_m(\alpha)$ and $t_m(u_\alpha)$ just differ by a logarithmic factor. Therefore, the multiple testing procedure behaves almost as well as the best test among the considered collection of tests.

(ii) The key point of the proof of Theorem 1 is an exponential inequality for $U$-statistics of order 2 due to Houdré and Reynaud-Bouret [11]. The same result could also be obtained with an inequality due to Giné, Latala and Zinn [10].

(iii) We prove in Section 5 that if $\mathcal{M}$ is finite, then for all $m = (l, D)$ in $\mathcal{M}$, $t_m(u_\alpha)$ is precisely of order $\sqrt{D \log(|\mathcal{M}|/\alpha)}/n$, where $|\mathcal{M}|$ denotes the cardinality of $\mathcal{M}$. This allows us to establish optimal uniform separation rates over various classes of alternatives. Furthermore, considering the problem of testing uniformity on $[0, 1]$, if we take a collection $\{S_m, m \in \mathcal{M}\}$ which only contains a finite number of spaces generated by constant piecewise functions, we can thus see that our procedure is very close to the one proposed by Ingster [15]. This would therefore be satisfactory enough from a theoretical point of view. Our choice to use a collection of mixing spaces generated by constant piecewise functions, scaling functions and possibly trigonometric polynomials is in fact explained by the experimental results. We indeed noticed in the simulation study that such a choice mostly increases the power of the test.

2.3. *Uniform separation rates.* Our purpose in this section is to evaluate the uniform separation rates of the test proposed above over several classes of alternatives. For $s > 0$, $R > 0$, $M > 0$ and $l \in \{1, 2, 3\}$, we introduce

$$\mathcal{B}_s^{(l)}(R, M)$$
$$= \{f \in \mathbb{L}_2(\mathbb{R}), \forall D \in \mathbb{D}_l, \|f - \Pi_{S_{(l,D)}}(f)\|_2^2 \le R^2 D^{-2s}, \|f\|_\infty \le M\}.$$

These sets of functions include some Hölder balls or Besov bodies. To be more precise, we consider, for all $s > 0$ and $R > 0$, the class of functions $\mathcal{H}_s(R)$ defined by

$$(2.4) \qquad \mathcal{H}_s(R) = \{f : [0, 1] \to \mathbb{R}, \forall x, y \in [0, 1],$$
$$|f^{(s_1)}(x) - f^{(s_1)}(y)| \le R|x - y|^{s_2}\},$$

where $s = s_1 + s_2$, $s_1 \in \mathbb{N}$ and $s_2 \in ]0, 1]$.

Let for $j$ in $\mathbb{N}$, $k$ in $\mathbb{Z}$, $\beta_{j,k}(f) = \langle f, \psi_{j,k} \rangle$. For all $s > 0$ and $R > 0$, we define the Besov body $B_{s,2,\infty}(R)$ as

$$B_{s,2,\infty}(R) = \left\{ f \in \mathbb{L}_2(\mathbb{R}), \ \forall j \in \mathbb{N}, \sum_{k \in \mathbb{Z}} \beta_{j,k}^2(f) \le R^2 2^{-2js} \right\}.$$

Then, one can see by straightforward computations that for $s \in ]0, 1]$, $R > 0$, $M > 0$,

$$\mathcal{H}_s(R) \cap \{f, \ \|f\|_\infty \le M\} \subset \mathcal{B}_s^{(1)}(R, M) \ \cap \ \mathcal{B}_s^{(3)}(R/\sqrt{2(4^s - 1)}, M),$$

and for $s > 0$, $R > 0$, $M > 0$,

$$B_{s,2,\infty}(R) \cap \{f, \; \|f\|_\infty \leq M\} \subset \mathcal{B}_s^{(2)}(R/\sqrt{1 - 4^{-s}}, M).$$

The following corollary gives upper bounds for the uniform separation rates of our testing procedure over the classes $\mathcal{B}_s^{(l)}(R, M)$.

COROLLARY 1. *Let $T_\alpha$ be the test statistic defined by* (2.3). *Assume that $n \geq 16$ and that for $l$ in $\{1, 2, 3\}$, $\mathcal{D}_l$ is $\{2^J, 0 \leq J \leq \log_2(n^2/(\log\log n)^3)\}$ or $\varnothing$. Fix some $\beta$ in $]0, 1[$. For all $s > 0$, $M > 0$, $R > 0$ and $l \in \{1, 2, 3\}$ such that $\mathcal{D}_l \neq \varnothing$, there exists some positive constant $C = C(s, \alpha, \beta, M, \|f_0\|_\infty)$ such that if $f$ belongs to the set $\mathcal{B}_s^{(l)}(R, M)$ and satisfies*

$$\|f - f_0\|_2^2 > C \bigg( R^{2/(4s+1)} \bigg( \frac{\sqrt{\log\log n}}{n} \bigg)^{4s/(4s+1)}$$
$$+ R^2 \bigg( \frac{(\log\log n)^3}{n^2} \bigg)^{2s} + \frac{(\log\log n)\log n}{n} \bigg),$$

*then*

$$\mathbb{P}_f(T_\alpha \leq 0) \leq \beta.$$

*In particular, if*

$$(2.5) \qquad (\log\log n)^{s+1/2}(\log n)^{2s+1/2}/\sqrt{n} \leq R \leq n^{2s}/(\log\log n)^{3s+1/2},$$

*there exists some positive constant $C'(s, \alpha, \beta, M, \|f_0\|_\infty)$ such that the uniform separation rate of the test $\mathbb{I}_{T_\alpha > 0}$ over $\mathcal{B}_s^{(l)}(R, M)$ satisfies*

$$\rho(\mathbb{I}_{T_\alpha > 0}, \mathcal{B}_s^{(l)}(R, M), \beta)$$
$$\leq C'(s, \alpha, \beta, M, \|f_0\|_\infty) R^{1/(4s+1)} \bigg( \frac{\sqrt{\log\log n}}{n} \bigg)^{2s/(4s+1)}.$$

COMMENTS. (i) For the problem of testing the null hypothesis "$f = \mathbb{I}_{[0,1]}$" against the alternative "$f = \mathbb{I}_{[0,1]} + g$ with $g \neq 0$ and $g \in B_s(R)$" where $B_s(R)$ is a class of smooth functions (like some Hölder, Sobolev or Besov ball in $\mathbb{L}_2([0,1])$) with unknown smoothness parameter $s$, Ingster [15] establishes that the adaptive minimax rate of testing is of order $(\sqrt{\log\log n}/n)^{2s/(4s+1)}$. From Corollary 1, we thus deduce that the procedure that we propose is adaptive in the sense that it is rate optimal over all the classes $\mathcal{B}_s^{(l)}(R, M)$ such that $R$ belongs to the range given by (2.5) simultaneously.

(ii) Ingster [15] considers in fact the minimax rates of testing with respect to general $\mathbb{L}_p$-distances. In particular, for $1 \leq p \leq 2$, he obtains the same adaptive minimax rate of testing $(\sqrt{\log\log n}/n)^{2s/(4s+1)}$. Our results can clearly be extended to $\mathbb{L}_p$-distances with $1 \leq p \leq 2$ when $f$ and $f_0$ have bounded support. In this case, one actually has that $\|f - f_0\|_p \leq C(p)\|f - f_0\|_2$.

We focus here on classes of alternatives that are well approximated by their projections onto the spaces $\{S_m, m \in \mathcal{M}\}$ under consideration. In the particular case where $f_0 = \mathbb{I}_{[0,1]}$, one can see that the test may be powerful even for alternatives that do not have such approximation properties. This is the purpose of Corollary 2.

COROLLARY 2.  *Let $f_0 = \mathbb{I}_{[0,1]}$. Assume that $n$ is larger than 16 and that $\mathcal{M} = \{(1, D), D \in \mathcal{D}_1\}$ with $\mathcal{D}_1 = \{2^J, 0 \leq J \leq \log_2(n^2/(\log \log n)^3)\}$. Let $T_\alpha$ be the test statistic defined by (2.3). For all $s > 0$ and $R > 0$, consider $\mathcal{H}_s(R)$ given by (2.4). Fix some $\beta$ in $]0, 1[$. For all $s > 0$, $M > 0$, $R > 0$, there exists some positive constant $C(R, s, \alpha, \beta, M)$ such that if $f$ belongs to the set $\mathcal{H}_s(R)$ with $\|f\|_\infty \leq M$, and if $f$ satisfies*

$$\|f - f_0\|_2^2 > C(R, s, \alpha, \beta, M)\left(\frac{\sqrt{\log \log n}}{n}\right)^{4s/(4s+1)},$$

*then*

$$\mathbb{P}_f(T_\alpha \leq 0) \leq \beta.$$

COMMENT.  The key point of the proof is an inequality due to Ingster ([14], part III, inequality (5.16)). This inequality allows in fact to avoid evaluating the approximation terms $\|f - \Pi_{S_m}(f)\|_2$. Although the functions $f$ in $\mathcal{H}_s(R)$ are not well approximated by their projections onto the spaces $\{S_{(1,D)}, D \in \mathcal{D}_1\}$ when $s > 1$, we thus prove that the corresponding testing procedure still achieves the adaptive minimax rate of testing over $\mathcal{H}_s(R)$ for any $s > 0$.

**3. Testing a parametric family.**  Let $X_1, \ldots, X_n$ be i.i.d. real random variables with common density $f$. In this section we consider the problem of testing that $f$ belongs to some translation/scale family of the form

$$\mathcal{F} = \left\{\frac{1}{\sigma}f_0\left(\frac{\cdot - \mu}{\sigma}\right), (\mu, \sigma) \in K\right\},$$

where $f_0$ is some given density and $K$ is some subset of $\mathbb{R} \times ]0, +\infty[$. The families of Gaussian, uniform or exponential densities and translation models are typical examples of such translation/scale families.

3.1. *Description of the test.*  The testing procedure introduced below is essentially based on the idea that if $f$ belongs to $\mathcal{F}$, there exists $(\mu, \sigma)$ in $K$ such that the density of the variables $(X_i - \mu)/\sigma$ is $f_0$. As in Section 2.1, we take an at most countable collection $\{S_m, m \in \mathcal{M}\}$ of linear subspaces of $\mathbb{L}_2(\mathbb{R})$. For all $m$ in $\mathcal{M}$, we consider an orthonormal basis $\{p_l, l \in \mathcal{L}_m\}$ of $S_m$ composed of right-continuous functions and we set

$$
\begin{aligned}
(3.1) \quad & \tilde{T}_m(X_1, \ldots, X_n) \\
& = \inf_{(\mu, \sigma) \in K} \Bigg\{ \frac{1}{n(n-1)} \sum_{l \in \mathcal{L}_m} \sum_{i \neq j=1}^{n} p_l\left(\frac{X_i - \mu}{\sigma}\right) p_l\left(\frac{X_j - \mu}{\sigma}\right) \\
& \hspace{4cm} + \|f_0\|_2^2 - \frac{2}{n} \sum_{i=1}^{n} f_0\left(\frac{X_i - \mu}{\sigma}\right) \Bigg\}.
\end{aligned}
$$

Since the functions $p_l$ are right-continuous, the infimum over $(\mu, \sigma)$ in $K$ can be replaced by the infimum over $(\mu, \sigma)$ in $K \cap \mathbb{Q}^2$ so that $\tilde{T}_m(X_1, \ldots, X_n)$ is a random variable.

We reject the null hypothesis "$f \in \mathcal{F}$" if

$$
\tilde{T}_\alpha = \sup_{m \in \mathcal{M}} (\tilde{T}_m(X_1, \ldots, X_n) - \tilde{q}_{m,\alpha})
$$

is positive, where $\{\tilde{q}_{m,\alpha}, \ m \in \mathcal{M}\}$ is a family of positive numbers such that

$$
(3.2) \qquad\qquad \sup_{f \in \mathcal{F}} \mathbb{P}_f(\tilde{T}_\alpha > 0) \leq \alpha.
$$

Let us explain how we choose $\{\tilde{q}_{m,\alpha}, \ m \in \mathcal{M}\}$. We distinguish two cases.

(i) The first one corresponds to the case where for all $m$ in $\mathcal{M}$, the variable $\tilde{T}_m(X_1, \ldots, X_n)$ defined by (3.1) satisfies

$$
(3.3) \quad \forall (\mu, \sigma) \in K \qquad \tilde{T}_m(X_1, \ldots, X_n) = \tilde{T}_m\left(\frac{X_1 - \mu}{\sigma}, \ldots, \frac{X_n - \mu}{\sigma}\right).
$$

This equality holds if, for instance, $K = \mathbb{R} \times ]0, +\infty[$, $K = \{0\} \times ]0, +\infty[$ or $K = \mathbb{R} \times \{1\}$. In this case, we take $\tilde{q}_{m,\alpha} = \tilde{t}_m(\tilde{u}_\alpha)$, where $\tilde{t}_m(u)$ is the $(1-u)$ quantile of $\tilde{T}_m(X_1, \ldots, X_n)$ under the hypothesis "$f = f_0$," and $\tilde{u}_\alpha$ is taken as

$$
\tilde{u}_\alpha = \sup \Bigg\{ u \in ]0, 1[, \ \mathbb{P}_{f_0}\Big( \sup_{m \in \mathcal{M}} (\tilde{T}_m(X_1, \ldots, X_n) - \tilde{t}_m(u)) > 0 \Big) \leq \alpha \Bigg\}.
$$

The quantities $\tilde{u}_\alpha$ and $\{\tilde{t}_m(\tilde{u}_\alpha), m \in \mathcal{M}\}$ can be estimated by Monte Carlo experiments. Let us see how this choice leads to inequality (3.2). Under the null hypothesis, there exists $(\mu, \sigma)$ in $K$ such that the density of the variables $(X_i - \mu)/\sigma$ is $f_0$. From (3.3), one can then deduce

$$
\begin{aligned}
& \mathbb{P}_f\Big( \sup_{m \in \mathcal{M}} (\tilde{T}_m(X_1, \ldots, X_n) - \tilde{t}_m(\tilde{u}_\alpha)) > 0 \Big) \\
& = \mathbb{P}_f\Big( \sup_{m \in \mathcal{M}} \Big( \tilde{T}_m\Big(\frac{X_1 - \mu}{\sigma}, \ldots, \frac{X_n - \mu}{\sigma}\Big) - \tilde{t}_m(\tilde{u}_\alpha) \Big) > 0 \Big) \\
& = \mathbb{P}_{f_0}\Big( \sup_{m \in \mathcal{M}} (\tilde{T}_m(X_1, \ldots, X_n) - \tilde{t}_m(\tilde{u}_\alpha)) > 0 \Big),
\end{aligned}
$$

and according to the definition of $\tilde{t}_m(\tilde{u}_\alpha)$, this probability is at most $\alpha$.

(ii) The second one corresponds to the case where (3.3) is not satisfied. This occurs, for instance, if $K$ is a compact set of $\mathbb{R} \times ]0, +\infty[$. Here, we take $\tilde{q}_{m,\alpha} = t_m(u_\alpha)$ where, as in Section 2.1, $t_m(u)$ is the $(1-u)$ quantile of the variable

$$\hat{T}_m(X_1, \ldots, X_n) = \frac{1}{n(n-1)} \sum_{l \in \mathcal{L}_m} \sum_{i \neq j=1}^{n} p_l(X_i) p_l(X_j) + \|f_0\|_2^2 - \frac{2}{n} \sum_{i=1}^{n} f_0(X_i)$$

under the assumption that the variables $X_1, \ldots, X_n$ are i.i.d. with common density $f_0$, and $u_\alpha$ is defined as

$$u_\alpha = \sup \left\{ u \in ]0,1[, \ \mathbb{P}_{f_0}\left( \sup_{m \in \mathcal{M}} (\hat{T}_m(X_1, \ldots, X_n) - t_m(u)) > 0 \right) \leq \alpha \right\}.$$

Inequality (3.2) also holds in this case: if $f$ belongs to $\mathcal{F}$, there exists $(\mu, \sigma)$ in $K$ such that $f = f_0((\cdot - \mu)/\sigma)/\sigma$. By definition of $\tilde{T}_m(X_1, \ldots, X_n)$, one has the inequality

$$\tilde{T}_m(X_1, \ldots, X_n) \leq \hat{T}_m\left( \frac{X_1 - \mu}{\sigma}, \ldots, \frac{X_n - \mu}{\sigma} \right).$$

Hence,

$$\mathbb{P}_f(\tilde{T}_\alpha > 0) \leq \mathbb{P}_f\left( \sup_{m \in \mathcal{M}} \left( \hat{T}_m\left( \frac{X_1 - \mu}{\sigma}, \ldots, \frac{X_n - \mu}{\sigma} \right) - t_m(u_\alpha) \right) > 0 \right).$$

Since the variables $(X_i - \mu)/\sigma$ have $f_0$ as common density, it follows from the definition of $u_\alpha$ that the above quantity is smaller than $\alpha$.

We shall remark that the choice of $\{\tilde{q}_{m,\alpha}, m \in \mathcal{M}\}$ proposed in (ii) may lead to a conservative procedure. It is therefore preferable to use the procedure proposed in (i) whenever condition (3.3) holds.

3.2. *The power of the test.* In the following, we use the same notation as in Section 2.2.

(i) For all $D$ in $\mathbb{D}_1 = \mathbb{N}^*$ and $m = (1, D)$, we define

$$\hat{T}_m(X_1, \ldots, X_n) = \frac{1}{n(n-1)} \sum_{k \in \mathbb{Z}} \sum_{i \neq j=1}^{n} I_{D,k}(X_i) I_{D,k}(X_j) + \|f_0\|_2^2 - \frac{2}{n} \sum_{i=1}^{n} f_0(X_i).$$

(ii) Choose the scaling function $\varphi$ such that it satisfies the Lipschitz condition

$$\forall x, y \in \mathbb{R} \qquad |\varphi(x) - \varphi(y)| \leq C_\varphi |x - y|.$$

For all $D = 2^J$ in $\mathbb{D}_2 = \{2^J, J \in \mathbb{N}\}$ and $m = (2, D)$, we define

$$\hat{T}_m(X_1, \ldots, X_n) = \frac{1}{n(n-1)} \sum_{k \in \mathbb{Z}} \sum_{i \neq j=1}^{n} \varphi_{J,k}(X_i) \varphi_{J,k}(X_j) + \|f_0\|_2^2 - \frac{2}{n} \sum_{i=1}^{n} f_0(X_i).$$

We recall that we do not specify the dependence on $\varphi$ in the involved constants.

Let $\mathcal{D}_1 \subset \mathbb{D}_1$ and $\mathcal{D}_2 \subset \mathbb{D}_2$ such that $\mathcal{D}_1 \cup \mathcal{D}_2 \neq \varnothing$ and let

$$\mathcal{M} = \{(l, D), l \in \{1, 2\}, D \in \mathcal{D}_l\}.$$

For all $m$ in $\mathcal{M}$, we set

$$\tilde{T}_m(X_1, \ldots, X_n) = \inf_{(\mu, \sigma) \in K} \hat{T}_m\left(\frac{X_1 - \mu}{\sigma}, \ldots, \frac{X_n - \mu}{\sigma}\right).$$

We consider the test statistic

$$(3.4) \qquad \tilde{T}_\alpha = \sup_{m \in \mathcal{M}} (\tilde{T}_m(X_1, \ldots, X_n) - \tilde{q}_{m,\alpha}),$$

where $\{\tilde{q}_{m,\alpha}, \ m \in \mathcal{M}\}$ is a family of positive numbers satisfying (3.2).

THEOREM 2. *Let $X_1, \ldots, X_n$ be i.i.d. real random variables with common density $f \in \mathbb{L}_\infty(\mathbb{R})$. Let*

$$\mathcal{F} = \left\{\frac{1}{\sigma} f_0\left(\frac{\cdot - \mu}{\sigma}\right), (\mu, \sigma) \in K\right\},$$

*where $f_0$ is some given bounded density and $K = [\underline{\mu}, \overline{\mu}] \times [\underline{\sigma}, \overline{\sigma}]$, $\underline{\mu}, \overline{\mu}, \underline{\sigma}, \overline{\sigma}$ being real numbers such that $\underline{\sigma} > 0$. Suppose that the following hypotheses hold:*

$(h_1)$ *There exists some constant $C_{f_0} > 0$ such that for all $x, y$ in the support of $f$, $(\mu, \sigma)$ in $K$, $(\mu', \sigma')$ in $K$,*

$$\left|f_0\left(\frac{x - \mu}{\sigma}\right) - f_0\left(\frac{y - \mu'}{\sigma'}\right)\right| \leq C_{f_0}\left|\frac{x - \mu}{\sigma} - \frac{y - \mu'}{\sigma'}\right|.$$

$(h_2)$ *There exist $\nu > 0$ and $c > 0$ such that for all $k \geq 2$,*

$$\mathbb{E}(|X_i|^k) \leq \frac{k!}{2}\nu c^{k-2}.$$

*Let $\tilde{T}_\alpha$ be the test statistic defined by (3.4). Assume that $\mathcal{D}_2 \neq \varnothing$ and that $n$ is large enough so that $n \geq 3$ and $n^2(\overline{\mu} - \underline{\mu}) \wedge (\overline{\sigma} - \underline{\sigma}) \geq 2$. Fix some $\beta$ in $]0, 1[$. For all $\varepsilon$ in $]0, 2[$, there exists some positive constant $C = C(\underline{\mu}, \overline{\mu}, \underline{\sigma}, \overline{\sigma}, \|f_0\|_\infty, C_{f_0}, \|f\|_\infty, \nu, c, \beta, \varepsilon)$ such that, setting*

$$\tilde{V}_D(\beta) = C\left(\frac{\sqrt{D}}{n}\sqrt{\log(n^2 D)} + \frac{D}{n^2}\log^2(n^2 D) + \frac{\log(n^2 D)}{n}\right),$$

*if $f$ satisfies*

$$\inf_{(\mu, \sigma) \in K} \sigma\left\|f - \frac{1}{\sigma}f_0\left(\frac{\cdot - \mu}{\sigma}\right)\right\|_2^2$$

$$\geq (1+\varepsilon) \inf_{D \in \mathcal{D}_2} \Bigg\{ \sup_{(\mu,\sigma) \in K} \| \sigma f(\sigma. + \mu) - \Pi_{S_{(2,D)}}(\sigma f(\sigma. + \mu)) \|_2^2$$

$$+ \tilde{V}_D(\beta) + \tilde{q}_{(2,D),\alpha} \Bigg\},$$

*then*

$$\mathbb{P}_f(\tilde{T}_\alpha \leq 0) \leq \beta.$$

COMMENTS. (i) In Theorem 2, we only consider the case where $K$ is compact. This is due to technical reasons: we have to evaluate the supremum over $(\mu, \sigma)$ in $K$ of $U$-statistics of order 2 depending on $(\mu, \sigma)$. By considering a finite grid on the compact set $K$, we reduce the problem to control of a finite number of these $U$-statistics and we can use the inequality due to Houdré and Reynaud-Bouret [11] again to control each of them.

(ii) The condition $(h_1)$ is satisfied by families of Gaussian densities (when $f_0(x) = e^{-x^2/2}/\sqrt{2\pi}$) whatever the support of $f$ and also by families of exponential densities (when $f_0(x) = e^{-x}\mathbb{I}_{x \geq 0}$ and $K \subset \{0\} \times ]0, +\infty[$) when the support of $f$ is included in $[0, +\infty[$. As for families of uniform densities (when $f_0 = \mathbb{I}_{[0,1]}$), the condition $(h_1)$ is not satisfied but the result still holds; to see this, we refer to a theorem stated in [9], where this condition is replaced by some $\mathbb{L}_2$-entropy with bracketing condition on $\mathcal{F}$.

(iii) As pointed out by a referee, the condition $(h_2)$ can be slightly weakened. Bernstein's inequality used in the proof of Theorem 2 actually still holds when $\mathbb{E}[e^{tX_i}] \leq e^{c't^2/2}$ for $0 \leq t \leq T$ (see [21], Section 2.2).

3.3. *Uniform separation rates.* As in Section 2.3, Theorem 2 allows us to evaluate the uniform separation rates of our test over classes of smooth functions. For all $s > 0$, $R > 0$, $M > 0$, we consider the set

$$\tilde{\mathcal{B}}_s(R, M) = \{f \in \mathbb{L}_2(\mathbb{R}), \|f\|_\infty \leq M, \ \forall D \in \mathbb{D}_2, \ \forall (\mu, \sigma) \in K,$$

$$\|\sigma f(\sigma. + \mu) - \Pi_{S_{(2,D)}}(\sigma f(\sigma. + \mu))\|_2^2 \leq R^2 \sigma^{1+2s} D^{-2s}\}.$$

Such a set contains, among others, the functions $f$ belonging to some Besov ball and satisfying the inequality $\|f\|_\infty \leq M$. To see this, in the notation of DeVore and Lorentz [7], we introduce, for all $h > 0$ and $r \in \mathbb{N}^*$,

$$\Delta_h^r(f, x) = \sum_{k=0}^r \binom{r}{k}(-1)^{r-k} f(x + kh).$$

The $r$th modulus of smoothness of $f$ in $\mathbb{L}_2(\mathbb{R})$ is defined by

$$\omega_r(f, t)_2 = \sup_{0 \leq h \leq t} \|\Delta_h^r(f, \cdot)\|_2.$$

Then $f$ belongs to the Besov ball $B_{s,2,\infty}(R)$ if for $r = [s] + 1$,

$$\sup_{t>0} t^{-s}\omega_r(f,t)_2 \leq R.$$

One can easily see that

$$\omega_r(\sigma f(\sigma. + \mu), t)_2 = \sigma^{1/2}\omega_r(f, \sigma t)_2.$$

Let us now recall an inequality due to DeVore, Jawerth and Popov [6]: if the wavelet $\psi$ satisfies that for all $j < r$, $\int x^j \psi(x)\,dx = 0$, then, for every function $g$ in $\mathbb{L}_2(\mathbb{R})$, for all $j \geq 0$,

$$\sum_{k\in\mathbb{Z}} \beta_{j,k}^2(g) \leq C\omega_r^2(g, 2^{-j})_2,$$

where $C$ is an absolute constant. Hence, if $f$ belongs to the Besov ball $B_{s,2,\infty}(R)$, for all $J \geq 0$, for all $(\mu, \sigma) \in K$,

$$\|\sigma f(\sigma. + \mu) - \Pi_{S_{(2,2^J)}}(\sigma f(\sigma. + \mu))\|_2^2 = \sum_{j\geq J}\sum_{k\in\mathbb{Z}}\beta_{j,k}^2(\sigma f(\sigma. + \mu))$$

$$\leq C(1-4^{-s})^{-1}R^2\sigma^{1+2s}2^{-2Js}.$$

If, in addition, $\|f\|_\infty \leq M$, then $f \in \tilde{\mathcal{B}}_s(C^{1/2}(1-4^{-s})^{-1/2}R, M)$.

The following corollary gives upper bounds for the uniform separation rates of the testing procedure over the classes $\tilde{\mathcal{B}}_s(R, M)$.

COROLLARY 3. *Assume that the conditions of Theorem 2 are satisfied. Let $\tilde{T}_\alpha$ be the test statistic defined by (3.4) with $\tilde{q}_{m,\alpha} = t_m(u_\alpha)$ as explained in Section 3.1 [case (ii)]. Choose $\mathcal{D}_1 \subset \{2^J, 0 \leq J \leq \log_2(n^2)\}$ and $\mathcal{D}_2 = \{2^J, 0 \leq J \leq \log_2(n^2/\log^3 n)\}$. Let $\beta \in ]0,1[$. For all $s > 0$, $M > 0$, $R > 0$, there exists a positive constant $C = C(\underline{\mu}, \overline{\mu}, \underline{\sigma}, \overline{\sigma}, \|f_0\|_\infty, C_{f_0}, M, \nu, c, \alpha, \beta, s)$ such that, if $f$ belongs to the set $\tilde{\mathcal{B}}_s(R, M)$ and satisfies*

$$\inf_{(\mu,\sigma)\in K}\left\|f - \frac{1}{\sigma}f_0\left(\frac{.-\mu}{\sigma}\right)\right\|_2^2$$

$$> C\left(R^{2/(4s+1)}\left(\frac{\sqrt{\log n}}{n}\right)^{4s/(4s+1)} + R^2\left(\frac{(\log n)^3}{n^2}\right)^{2s} + \frac{\log n}{n}\right),$$

*then*

$$\mathbb{P}_f(\tilde{T}_\alpha \leq 0) \leq \beta.$$

COMMENTS. (i) When $R$ satisfies $(\log n)^{s+1/2}/\sqrt{n} \leq R \leq n^{2s}/(\log n)^{3s+1/2}$, the uniform separation rate over the class of functions belonging to $\tilde{\mathcal{B}}_s(R, M)$ and satisfying $(h_1)$ and $(h_2)$ is bounded from above by

$$C'(\underline{\mu}, \overline{\mu}, \underline{\sigma}, \overline{\sigma}, \|f_0\|_\infty, C_{f_0}M, \nu, c, \alpha, \beta, s)R^{1/(4s+1)}\left(\frac{\sqrt{\log n}}{n}\right)^{2s/(1+4s)}.$$

This corresponds, up to a logarithmic factor, to the rate over the classes $\mathcal{B}_s^{(l)}(R, M)$ for the test of simple hypotheses obtained in Corollary 1. We do not know if this logarithmic factor can be avoided.

(ii) As in Corollary 1, the result can be extended to $\mathbb{L}_p$-distances with $p$ in $[1, 2]$ when $f$ and $f_0$ have bounded support.

## 4. Simulation study.

4.1. *Test of uniformity on* $[0, 1]$. We first present simulation results for the problem of testing that the distribution of some i.i.d. random variables $X_1, \ldots, X_n$ with values in $[0, 1]$ is uniform on $[0, 1]$. In order to implement our procedure, we have to choose the set $\mathcal{M} = \{(l, D), l \in \{1, 2, 3\}, D \in \mathcal{D}_l\}$ that occurs in the definition (2.3) of the test statistic $T_\alpha$. We present two cases. In the first case, we consider only trigonometric polynomials. We take $\mathcal{D}_1 = \mathcal{D}_2 = \varnothing$ and $\mathcal{D}_3 = \{1, 2, \ldots, D_{\mathrm{tr}}\}$. Setting $g_0 = \mathbb{I}_{[0,1]}$ and for all $p \geq 1$, $g_{2p-1}(x) = \sqrt{2} \cos(2p\pi x) \mathbb{I}_{[0,1]}(x)$, and $g_{2p}(x) = \sqrt{2} \sin(2p\pi x) \mathbb{I}_{[0,1]}(x)$, the test statistic is based on orthogonal projections onto the spaces spanned by the functions $\{g_l, l = 0, \ldots, D\}$ for $D$ in $\mathcal{D}_3$. The second case consists in mixing trigonometric polynomials and constant piecewise functions: we take $\mathcal{D}_1 = \{2, 3, \ldots, D_{\mathrm{ct}}\}$ and $\mathcal{D}_3 = \{1, \ldots, D_{\mathrm{tr}}\}$. The tests corresponding to these two cases are, respectively, denoted by $\mathcal{T}_{\mathrm{tr}}$ and $\mathcal{T}_{\mathrm{tr/ct}}$. We compare their powers with the powers of the tests proposed by Kallenberg and Ledwina [17] (denoted by $\mathcal{T}_{\mathrm{KL}}$), Bickel and Ritov [4] (denoted by $\mathcal{T}_{\mathrm{BR}}$) and Kolmogorov and Smirnov (denoted by $\mathcal{T}_{\mathrm{KS}}$). As explained in the Introduction, Kallenberg and Ledwina propose a test of uniformity on $[0, 1]$ which is a data-driven version of Neyman's test [20]. They consider the orthonormal system $\{\phi_l, \ l \geq 0\}$ of $\mathbb{L}_2([0, 1])$, where $\phi_0 = \mathbb{I}_{[0,1]}$ and the $\phi_l$'s for $l \geq 1$ are the orthonormal Legendre polynomials on $[0, 1]$. They decide to reject the null hypothesis if the statistic $T_D = \sum_{l=1}^D (n^{-1/2} \sum_{i=1}^n \phi_l(X_i))^2$ is large, $D$ being chosen in $\{1, \ldots, d(n)\}$ via Schwarz's Bayesian Information Criterion. The critical value is estimated by simulations. The test proposed by Bickel and Ritov is based on the statistic

$$T_n = \max_{1 \leq D \leq d(n)} (T_{n,D} - D)/\sqrt{2D},$$

where $T_{n,D} = \frac{1}{n} \sum_{l=1}^D \sum_{i,j=1}^n 2 \cos(l\pi X_i) \cos(l\pi X_j)$.

We focus on the alternatives studied in the paper by Kallenberg and Ledwina [17], for which the power of Bickel and Ritov's test is also given. These alternatives are

$$f_{(\rho, j)}(x) = 1 + \rho \cos(j\pi x),$$
$$g_{(p, q, \varepsilon)}(x) = 1 - \varepsilon + \varepsilon \beta_{p,q}(x),$$
$$h_{(\rho, j)}(x) = 1 + \rho \phi_j(x),$$

where $\beta_{p,q}$ is the Beta density with parameter $(p,q)$ and $\{\phi_j, j \geq 1\}$ is the family of orthonormal Legendre polynomials on $[0,1]$.

We have chosen a level $\alpha = 5\%$.

The value of $u_\alpha$ and the quantiles $\{t_m(u_\alpha), m \in \mathcal{M}\}$ are estimated by 40,000 simulations. We use 20,000 simulations for the estimation of the $(1-u)$ quantiles $t_m(u)$ of the variables $\hat{T}_m = \hat{\theta}_m + \|f_0\|_2^2 - 2n^{-1}\sum_{i=1}^{n} f_0(X_i)$ under the hypothesis "$f = f_0$" for $u$ varying on a regular grid of $]0, \alpha[$ and 20,000 simulations for the estimation of the probabilities $\mathbb{P}_{f_0}(\sup_{m \in \mathcal{M}}(\hat{T}_m - t_m(u)) > 0)$.

Tables 1 and 2 present the estimated powers for the tests $\mathcal{T}_{\mathrm{tr}}$, $\mathcal{T}_{\mathrm{tr/ct}}$, $\mathcal{T}_{\mathrm{KL}}$, $\mathcal{T}_{\mathrm{BR}}$ and $\mathcal{T}_{\mathrm{KS}}$ under various alternatives for a number of observations equal to 50 or 100. The powers of the tests $\mathcal{T}_{\mathrm{tr}}$, $\mathcal{T}_{\mathrm{tr/ct}}$ and $\mathcal{T}_{\mathrm{KS}}$ are estimated by 5000 experiments and the levels by 20,000 experiments. Hence, with confidence 95%, the estimation error is less than 0.3% for the levels and less than 1.3% for the powers.

For a number of observations equal to 50, Kallenberg and Ledwina take $d(50) = 10$. We choose $D_{\mathrm{tr}} = 6$, $D_{\mathrm{ct}} = 6$, and we obtain the results in Table 1.

TABLE 1
*Estimated powers of the test of uniformity on $[0,1]$ for $n = 50$ with $D_{\mathrm{tr}} = 6$ and $D_{\mathrm{ct}} = 6$*

| **Alternatives $f_{(\rho,j)}$** | | | | |
|---|---|---|---|---|
| $(\rho,j)$ | $\mathcal{T}_{\mathrm{tr}}$ | $\mathcal{T}_{\mathrm{tr/ct}}$ | $\mathcal{T}_{\mathrm{KL}}$ | $\mathcal{T}_{\mathrm{BR}}$ | $\mathcal{T}_{\mathrm{KS}}$ |
| $(0.5, 2)$ | 0.61 | 0.56 | 0.56 | 0.48 | 0.29 |
| $(0.7, 4)$ | 0.80 | 0.77 | 0.50 | 0.71 | 0.16 |
| $(0.7, 6)$ | 0.69 | 0.62 | 0.23 | 0.60 | 0.10 |
| **Alternatives $g_{(p,q,\varepsilon)}$** | | | | |
| $(p,q,\varepsilon)$ | $\mathcal{T}_{\mathrm{tr}}$ | $\mathcal{T}_{\mathrm{tr/ct}}$ | $\mathcal{T}_{\mathrm{KL}}$ | $\mathcal{T}_{\mathrm{BR}}$ | $\mathcal{T}_{\mathrm{KS}}$ |
| $(3, 3, 1/2)$ | 0.55 | 0.49 | 0.53 | 0.40 | 0.14 |
| $(10, 20, 0.25)$ | 0.46 | 0.49 | 0.36 | 0.41 | 0.33 |
| $(2, 2, 0.8)$ | 0.62 | 0.55 | 0.63 | 0.44 | 0.15 |
| $(2, 4, 0.5)$ | 0.57 | 0.60 | 0.55 | 0.58 | 0.64 |
| **Alternatives $h_{(\rho,j)}$** | | | | |
| $(\rho,j)$ | $\mathcal{T}_{\mathrm{tr}}$ | $\mathcal{T}_{\mathrm{tr/ct}}$ | $\mathcal{T}_{\mathrm{KL}}$ | $\mathcal{T}_{\mathrm{BR}}$ | $\mathcal{T}_{\mathrm{KS}}$ |
| $(0.4, 2)$ | 0.69 | 0.65 | 0.70 | 0.59 | 0.32 |
| $(0.3, 5)$ | 0.16 | 0.16 | 0.13 | 0.14 | 0.07 |
| **Estimated levels** | | | | |
| $\mathcal{T}_{\mathrm{tr}}$ | $\mathcal{T}_{\mathrm{tr/ct}}$ | $\mathcal{T}_{\mathrm{KL}}$ | $\mathcal{T}_{\mathrm{BR}}$ | $\mathcal{T}_{\mathrm{KS}}$ |
| 0.051 | 0.055 | 0.061 | 0.031 | 0.050 |

TABLE 2
*Estimated powers of the test of uniformity on $[0,1]$ for $n = 100$ with $D_{tr} = 12$ and $D_{ct} = 10$*

| Alternatives $f_{(\rho,j)}$ | | | | |
|---|---|---|---|---|
| $(\rho, j)$ | $\mathcal{T}_{\text{tr}}$ | $\mathcal{T}_{\text{tr/ct}}$ | $\mathcal{T}_{\text{KL}}$ | $\mathcal{T}_{\text{BR}}$ | $\mathcal{T}_{\text{KS}}$ |
| $(0.5, 2)$ | 0.87 | 0.85 | 0.87 | 0.84 | 0.53 |
| $(0.7, 4)$ | 0.98 | 0.98 | 0.83 | 0.98 | 0.29 |
| $(0.7, 6)$ | 0.97 | 0.96 | 0.46 | 0.95 | 0.19 |
| **Alternatives $g_{(p,q,\varepsilon)}$** | | | | |
| $(p, q, \varepsilon)$ | $\mathcal{T}_{\text{tr}}$ | $\mathcal{T}_{\text{tr/ct}}$ | $\mathcal{T}_{\text{KL}}$ | $\mathcal{T}_{\text{BR}}$ | $\mathcal{T}_{\text{KS}}$ |
| $(3, 3, 1/2)$ | 0.83 | 0.77 | 0.88 | 0.76 | 0.35 |
| $(10, 20, 0.25)$ | 0.77 | 0.78 | 0.62 | 0.75 | 0.60 |
| $(2, 2, 0.8)$ | 0.90 | 0.86 | 0.95 | 0.82 | 0.36 |
| $(2, 4, 0.5)$ | 0.87 | 0.89 | 0.88 | 0.90 | 0.91 |
| **Alternatives $h_{(\rho,j)}$** | | | | |
| $(\rho, j)$ | $\mathcal{T}_{\text{tr}}$ | $\mathcal{T}_{\text{tr/ct}}$ | $\mathcal{T}_{\text{KL}}$ | $\mathcal{T}_{\text{BR}}$ | $\mathcal{T}_{\text{KS}}$ |
| $(0.4, 2)$ | 0.93 | 0.91 | 0.95 | 0.90 | 0.60 |
| $(0.3, 5)$ | 0.33 | 0.31 | 0.23 | 0.33 | 0.09 |
| **Estimated levels** | | | | |
| $\mathcal{T}_{\text{tr}}$ | $\mathcal{T}_{\text{tr/ct}}$ | $\mathcal{T}_{\text{KL}}$ | $\mathcal{T}_{\text{BR}}$ | $\mathcal{T}_{\text{KS}}$ |
| 0.050 | 0.048 | 0.056 | 0.031 | 0.054 |

For a number of observations equal to 100, Kallenberg and Ledwina take $d(100) = 12$. We choose $D_{tr} = 12$, $D_{ct} = 10$, and we obtain the results in Table 2.

COMMENTS. In this simulation study, the alternatives that we consider are of three kinds. The alternatives $f_{(\rho,j)}$ correspond to the uniform density contaminated by a cosine function. They are favorable to our tests and Bickel and Ritov's test since these tests are based on trigonometric polynomials. It is therefore natural to compare our power results with the results of $\mathcal{T}_{\text{BR}}$. The main difference between the two procedures lies in the fact that $\hat{\theta}_{(3,D)}$ is an unbiased estimator of the squared $\mathbb{L}_2$-norm of the orthogonal projection of $f$ onto $S_{(3,D)}$, whereas $(T_{n,D} - D)/(n-1)$ is an unbiased estimator of $\sum_{l=1}^{D} (\int f(x)\sqrt{2}\cos(l\pi x)\,dx)^2$ only under the null hypothesis. The consequent bias term under some alternative $f$ may be of order $D/n$, which does not allow us to establish optimal uniform separation rates for the test proposed by Bickel and Ritov. This explains why the power of $\mathcal{T}_{\text{tr}}$ improves the power of $\mathcal{T}_{\text{BR}}$ in all cases.

The alternatives $h_{(\rho,j)}$ are more favorable to the test due to Kallenberg and Ledwina since this test is based on Legendre polynomials. However, under the alternative $h_{(0.3,5)}$, the test $\mathcal{T}_{tr}$ improves the results of $\mathcal{T}_{KL}$ and under the alternative $h_{(0.4,2)}$, the estimated power of $\mathcal{T}_{tr}$ is comparable to that of $\mathcal{T}_{KL}$.

Since the functions $g_{(p,q,\varepsilon)}$ correspond to some "neutral" alternatives, we can focus on them. When the number of observations is equal to 100, the powers of our tests $\mathcal{T}_{tr}$ and $\mathcal{T}_{tr/ct}$ are at least equivalent to the powers of $\mathcal{T}_{KL}$ and $\mathcal{T}_{BR}$ for half of the considered cases. As for the other cases, the procedures $\mathcal{T}_{tr}$ and $\mathcal{T}_{tr/ct}$ are still more powerful than $\mathcal{T}_{BR}$. For small sample sizes ($n = 50$), our test $\mathcal{T}_{tr}$ is always at least as powerful as the tests $\mathcal{T}_{KL}$ and $\mathcal{T}_{BR}$.

4.2. *Other goodness-of-fit tests.* We are now interested in testing that the density $f$ of the random variables $X_1, \ldots, X_n$ is a given density $f_0$, with $f_0 \neq \mathbb{I}_{[0,1]}$. To test such a hypothesis, we have two possible procedures: the first one consists in testing directly from the sample $X_1, \ldots, X_n$ the null hypothesis "$f = f_0$" as explained in Section 2.1. The second one consists in testing that the common distribution of the variables $F_0(X_1), \ldots, F_0(X_n)$, where $F_0$ is the distribution function associated with the density $f_0$, is uniform on $[0,1]$. This approach is the one which is proposed in most papers. Whereas the two procedures are equivalent for the Kolmogorov–Smirnov test, they are not for our method based on the estimation of an $\mathbb{L}_2$-distance. Indeed, in our case, the first test is based on the estimation of $\|f - f_0\|_2^2$. Since the density of $F_0(X_1)$ is given by $h(x) = f(F_0^{-1}(x))/f_0(F_0^{-1}(x))$ when $F_0$ is one to one, the second test is based on the estimation of

$$\int_{[0,1]} (h(x) - 1)^2 \, dx = \|(f - f_0)/\sqrt{f_0}\|_2^2.$$

To compare in practice these two procedures, we have chosen to test that the density $f$ is Gaussian with mean 0 and with variance 1 first, with variance 0.01 second.

The choices we have made in order to implement our procedures are the following ones: for the direct test from $X_1, \ldots, X_n$ denoted by $\mathcal{T}_d$, the set $\mathcal{M} = \{(l,D), l \in \{1,2,3\}, D \in \mathcal{D}_l\}$ is taken such that $\mathcal{D}_2 = \mathcal{D}_3 = \varnothing$ and $\mathcal{D}_1 = \{1, \ldots, 10\}$, and for the second test from $F_0(X_1), \ldots, F_0(X_n)$, we apply the test $\mathcal{T}_{tr/ct}$ described in Section 4.1 with $D_{ct} = 10$ and $D_{tr} = 12$. We also present the estimated powers of the Kolmogorov–Smirnov test that we still denote by $\mathcal{T}_{KS}$.

We have taken a number of observations $n = 100$ and a level $\alpha = 5\%$. The quantiles are estimated as above with 40,000 simulations, the powers of the tests with 5000 simulations and the levels with 20,000 simulations. The

TABLE 3
*Estimated powers for $n = 100$*

| Test of normality $\mathcal{N}(0, 1)$ | | | | Test of normality $\mathcal{N}(0, 0.01)$ | | | |
|---|---|---|---|---|---|---|---|
| **Alternatives $f_m$** | | | | | | | |
| $m$ | $\mathcal{T}_d$ | $\mathcal{T}_{tr/ct}$ | $\mathcal{T}_{KS}$ | $m$ | $\mathcal{T}_d$ | $\mathcal{T}_{tr/ct}$ | $\mathcal{T}_{KS}$ |
| 2 | 0.96 | 0.92 | 0.62 | 0.17 | 0.93 | 0.64 | 0.24 |
| 1.8 | 0.66 | 0.66 | 0.36 | 0.16 | 0.87 | 0.71 | 0.14 |
| $\sqrt{\pi/2}$ | 0.71 | 1 | 0.07 | 0.12 | 0.99 | 1 | 0.14 |
| **Alternatives $g_{(m,\sigma^2)}$** | | | | | | | |
| $(m, \sigma^2)$ | $\mathcal{T}_d$ | $\mathcal{T}_{tr/ct}$ | $\mathcal{T}_{KS}$ | $(m, \sigma^2)$ | $\mathcal{T}_d$ | $\mathcal{T}_{tr/ct}$ | $\mathcal{T}_{KS}$ |
| $(1, 1)$ | 0.80 | 0.98 | 0.77 | $(0.1, 0.01)$ | 1 | 0.98 | 0.77 |
| $(0.5, 2)$ | 0.66 | 0.98 | 0.70 | $(0.05, 0.015)$ | 0.91 | 0.77 | 0.35 |
| $(1, 2)$ | 0.97 | 1 | 0.97 | $(0.05, 0.02)$ | 1 | 0.97 | 0.68 |
| **Alternatives $h_p$** | | | | | | | |
| $p$ | $\mathcal{T}_d$ | $\mathcal{T}_{tr/ct}$ | $\mathcal{T}_{KS}$ | $p$ | $\mathcal{T}_d$ | $\mathcal{T}_{tr/ct}$ | $\mathcal{T}_{KS}$ |
| $2/\sqrt{2\pi}$ | 0.24 | 0.95 | 0.42 | $20/\sqrt{2\pi}$ | 0.96 | 0.95 | 0.41 |
| $3/2\sqrt{2\pi}$ | 0.85 | 1 | 0.96 | $15/\sqrt{2\pi}$ | 1 | 1 | 0.96 |
| **Estimated levels** | | | | | | | |
| | $\mathcal{T}_d$ | $\mathcal{T}_{tr/ct}$ | $\mathcal{T}_{KS}$ | | $\mathcal{T}_d$ | $\mathcal{T}_{tr/ct}$ | $\mathcal{T}_{KS}$ |
| | 0.052 | 0.051 | 0.053 | | 0.053 | 0.055 | 0.053 |

alternatives that we have considered are the following ones (see Table 3):

$$f_m(x) = \frac{1}{2m} \mathbb{I}_{[-m,m]},$$

$$g_{(m,\sigma^2)}(x) = \frac{1}{2\sqrt{2\pi}\sigma}(e^{-(x-m)^2/(2\sigma^2)} + e^{-(x+m)^2/(2\sigma^2)}),$$

$$h_p(x) = \frac{p}{2}e^{px}\mathbb{I}_{x<0} + \frac{p}{2}e^{-px}\mathbb{I}_{x\geq 0}.$$

COMMENTS. The first objective of this simulation study is to compare our tests with the Kolmogorov–Smirnov test. The estimated power of the most powerful of our tests is larger than that of $\mathcal{T}_{KS}$ for most of the considered alternatives. In these cases, the difference in power is really significant as we can see, for example, for the alternatives $f_m$.

The second objective is to compare the tests $\mathcal{T}_d$ and $\mathcal{T}_{tr/ct}$. We can notice that we reject the null hypothesis more often with the test of uniformity from the $F_0(X_i)$'s when we test that the density is Gaussian with variance 1, whereas the direct test performs better when we test that the density is Gaussian with variance 0.01. As explained above, this is due to the fact that $\|(f - f_0)/\sqrt{f_0}\|_2^2$ is larger than $\|f - f_0\|_2^2$ when $f_0$ is the standard Gaussian density but smaller when $f_0$ is the Gaussian density with variance 0.01.

4.3. *Testing a parametric family.* We now implement the testing procedure described in Section 3 in order to test that the density $f$ of the observations $X_1, \ldots, X_n$ is an exponential density or, in other words, that $f$ belongs to the set of densities

$$\mathcal{F} = \{f, \ f(x) = \sigma^{-1} e^{-x/\sigma} \mathbb{I}_{x \geq 0}, \ \sigma > 0\}.$$

To simplify the implementation, we base our test statistic on constant piecewise functions instead of scaling functions. For all $D$ in $\{2, \ldots, 10\}$, we define

$$\tilde{T}_{(1,D)}(X_1, \ldots, X_n)$$
$$= \inf_{\sigma > 0} \left\{ \frac{D}{n(n-1)} \sum_{k \geq 0} \sum_{i \neq j = 1}^{n} \mathbb{I}_{\{X_i/\sigma \in [k/D, (k+1)/D[\}} \mathbb{I}_{\{X_j/\sigma \in [k/D, (k+1)/D[\}} \right.$$
$$\left. + \|f_0\|_2^2 - \frac{2}{n} \sum_{i=1}^{n} f_0\left(\frac{X_i}{\sigma}\right) \right\},$$

where $f_0$ is given by $f_0(x) = e^{-x} \mathbb{I}_{x \geq 0}$. It is easy to see that for all $\sigma > 0$, $\tilde{T}_{(1,D)}(X_1/\sigma, \ldots, X_n/\sigma) = \tilde{T}_{(1,D)}(X_1, \ldots, X_n)$.

Let $\mathcal{M} = \{(1, D), D = 2, \ldots, 10\}$. As explained in Section 3, we obtain a level $\alpha$ test as follows: denoting by $\tilde{t}_m(u)$ the $(1-u)$ quantile of $\tilde{T}_m(X_1, \ldots, X_n)$ under the hypothesis "$f = f_0$," and setting

$$\tilde{u}_\alpha = \sup \left\{ u \in ]0, 1[, \mathbb{P}_{f_0}\left( \sup_{m \in \mathcal{M}} (\tilde{T}_m(X_1, \ldots, X_n) - \tilde{t}_m(u)) > 0 \right) \leq \alpha \right\},$$

we reject the null hypothesis if

$$\sup_{m \in \mathcal{M}} (\tilde{T}_m(X_1, \ldots, X_n) - \tilde{t}_m(\tilde{u}_\alpha)) > 0.$$

The quantities $\tilde{u}_\alpha$ and $\{\tilde{t}_m(\tilde{u}_\alpha), m \in \mathcal{M}\}$ are estimated by Monte Carlo experiments: 20,000 simulations are used to estimate the values of $\{\tilde{t}_m(u), m \in \mathcal{M}\}$ for $u$ varying on a regular grid of $]0, \alpha[$, and 20,000 simulations are used to estimate the value of $\tilde{u}_\alpha$.

We compare the performance of our test with that of the Kolmogorov–Smirnov test described below. Let

$$\bar{X}_n = n^{-1} \sum_{i=1}^{n} X_i$$

and let $\hat{F}_n$ be the empirical distribution function defined for all $t$ in $\mathbb{R}$ by

$$\hat{F}_n(t) = n^{-1} \sum_{i=1}^{n} \mathbb{I}_{X_i \leq t}.$$

Let $F_\sigma$ be the distribution function associated with the exponential density with parameter $\sigma^{-1} : F_\sigma(t) = (1 - e^{-t/\sigma})\mathbb{I}_{t \geq 0}$. Under the hypothesis that the distribution of the $X_i$'s is exponential with parameter $\sigma^{-1}$, the law of the statistic $D_n = \sup_{t \in \mathbb{R}} |\hat{F}_n(t) - F_{\bar{X}_n}(t)|$ is free of the parameter $\sigma$. Let $d_{n,1-\alpha}$ denote the $(1-\alpha)$ quantile of $D_n$ under the assumption that the $X_i$'s have $f_0$ as common density. The Kolmogorov–Smirnov test of exponentiality consists in rejecting the null hypothesis "$f \in \mathcal{F}$" if $D_n > d_{n,1-\alpha}$. We consider the alternatives defined for $x > 0$ by

$$g_p(x) = (e^{-x} + (1 + \sin(p\pi x))\mathbb{I}_{0 < x < 1})/2 \qquad (p \text{ even}),$$

$$h_p(x) = (e^{-x} + (1 + \cos(p\pi x))\mathbb{I}_{0 < x < 1})/2,$$

$$k_{(p,q,\varepsilon)}(x) = (1 - \varepsilon)e^{-x} + \varepsilon\beta_{p,q}(x),$$

$$l_{(p,q,\varepsilon)}(x) = (1 - \varepsilon)e^{-x} + \varepsilon\gamma_{p,q}(x),$$

$$t(x) = e^{-(\log x)^2/2}/(x\sqrt{2\pi}),$$

$$v(x) = \sqrt{x}e^{-x/2}/(2^{3/2}\Gamma(3/2)),$$

$$w(x) = 1.5x^{0.5}e^{-x^{1.5}},$$

where $\beta_{p,q}$ and $\gamma_{p,q}$, respectively, denote the Beta density and the Gamma density with parameters with parameters $(p, q)$.

For each alternative, the power of the test is still estimated by 5000 simulations. The levels are estimated by 20,000 simulations. We choose $n = 100$ and $\alpha = 5\%$.

Table 4 presents the estimated power of our test denoted by $\mathcal{T}'_\alpha$ and of the Kolmogorov–Smirnov test denoted by $\mathcal{T}'_{\text{KS}}$.

COMMENT. We can see in Table 4 that our test is not always much more powerful than the Kolmogorov–Smirnov test under very smooth alternatives like $v$ and $w$ which respectively correspond to the chi-square with three degrees of freedom and the Weibull with parameter 1.5. However, under oscillating alternatives like $g_p$, $h_p$, $k_{(p,q,\varepsilon)}$ and $l_{(p,q,\varepsilon)}$, for which the Kolmogorov–Smirnov test is known to fail, our test performs much better. We could furthermore expect better results for our procedure with regular scaling functions instead of constant piecewise functions that we have chosen to make the implementation easier.

**5. Proofs.**

TABLE 4
*Estimated power of the test of exponentiality for $n = 100$*

| | Alternatives $g_p$ | | Alternative $t$ | |
|---|---|---|---|---|
| $p$ | $\mathcal{T}'_\alpha$ | $\mathcal{T}'_{\mathrm{KS}}$ | $\mathcal{T}'_\alpha$ | $\mathcal{T}'_{\mathrm{KS}}$ |
| 4 | 0.89 | 0.74 | 0.75 | 0.45 |
| | Alternatives $h_p$ | | Alternative $v$ | |
| $p$ | $\mathcal{T}'_\alpha$ | $\mathcal{T}'_{\mathrm{KS}}$ | $\mathcal{T}'_\alpha$ | $\mathcal{T}'_{\mathrm{KS}}$ |
| 4 | 0.71 | 0.60 | 0.67 | 0.65 |
| 1 | 1 | 0.90 | | |
| | Alternatives $k_{(p,q,\varepsilon)}$ | | Alternative $w$ | |
| $(p, q, \varepsilon)$ | $\mathcal{T}'_\alpha$ | $\mathcal{T}'_{\mathrm{KS}}$ | $\mathcal{T}'_\alpha$ | $\mathcal{T}'_{\mathrm{KS}}$ |
| (10, 20, 0.25) | 0.91 | 0.65 | 0.97 | 0.98 |
| | Alternatives $l_{(p,q,\varepsilon)}$ | | Estimated levels | |
| $(p, q, \varepsilon)$ | $\mathcal{T}'_\alpha$ | $\mathcal{T}'_{\mathrm{KS}}$ | $\mathcal{T}'_\alpha$ | $\mathcal{T}'_{\mathrm{KS}}$ |
| (2, 5, 0.5) | 0.53 | 0.28 | 0.053 | 0.051 |
| (2, 5, 0.75) | 0.89 | 0.60 | | |

5.1. *Proof of Theorem* 1. The main tool of the proof is the canonical decomposition of the $U$-statistics $\hat{\theta}_m$ defined in Section 2.2. We introduce the processes $U_n$ and $P_n$ defined by

$$U_n(h) = \frac{1}{n(n-1)} \sum_{i \neq j=1}^{n} h(X_i, X_j), \qquad P_n(h) = \frac{1}{n} \sum_{i=1}^{n} h(X_i).$$

We also define $P(h) = \langle h, f \rangle$. Using the same notation as in Section 2.2, let for $D$ in $\mathbb{N}^*$, $J$ in $\mathbb{N}$,

(5.1)
(i) $\mathcal{L}_{(1,D)} = \{D\} \times \mathbb{Z}$ and $\{p_l, l \in \mathcal{L}_{(1,D)}\} = \{I_{D,k}, k \in \mathbb{Z}\}$,
(ii) $\mathcal{L}_{(2,2^J)} = \{J\} \times \mathbb{Z}$ and $\{p_l, l \in \mathcal{L}_{(2,2^J)}\} = \{\varphi_{J,k}, k \in \mathbb{Z}\}$,
(iii) $\mathcal{L}_{(3,D)} = \{0, 1, \ldots, D\}$ and $\{p_l, l \in \mathcal{L}_{(3,D)}\} = \{g_l, l = 0, \ldots, D\}$.

By setting, for all $m \in \mathcal{M}$,

$$H_m(x,y) = \sum_{l \in \mathcal{L}_m} (p_l(x) - a_l)(p_l(y) - a_l),$$

with $a_l = \langle f, p_l \rangle$, we obtain the decomposition

$$\hat{\theta}_m = U_n(H_m) + (P_n - P)(2\Pi_{S_m}(f)) + \|\Pi_{S_m}(f)\|_2^2.$$

Let us fix some $\beta$ in $]0, 1[$. Recalling that

$$\mathbb{P}_f(T_\alpha \leq 0) = \mathbb{P}_f\left(\sup_{m \in \mathcal{M}}\left(\hat{\theta}_m + \|f_0\|_2^2 - \frac{2}{n}\sum_{i=1}^n f_0(X_i) - t_m(u_\alpha)\right) \leq 0\right),$$

we have

$$(5.2) \quad \mathbb{P}_f(T_\alpha \leq 0) \leq \inf_{m \in \mathcal{M}} \mathbb{P}_f\left(\hat{\theta}_m + \|f_0\|_2^2 - \frac{2}{n}\sum_{i=1}^n f_0(X_i) - t_m(u_\alpha) \leq 0\right).$$

Since $\|f - \Pi_{S_m}(f)\|_2^2 = \|f\|_2^2 - \|\Pi_{S_m}(f)\|_2^2$,

$$\begin{aligned}
(5.3) \quad \mathbb{P}_f(T_\alpha \leq 0) \leq \inf_{m \in \mathcal{M}} \mathbb{P}_f(&U_n(H_m) + (P_n - P)(2\Pi_{S_m}(f) - 2f) \\
&+ (P_n - P)(2f - 2f_0) \\
&+ \|f - f_0\|_2^2 \leq \|f - \Pi_{S_m}(f)\|_2^2 + t_m(u_\alpha)).
\end{aligned}$$

We then need to control $U_n(H_m)$, $(P_n - P)(2\Pi_{S_m}(f) - 2f)$, $(P_n - P)(2f - 2f_0)$ for every $m$ in $\mathcal{M}$.

(a) Control of $U_n(H_m)$.   We use the following lemma, which derives from an exponential inequality for $U$-statistics of order 2 due to Houdré and Reynaud-Bouret [11].

LEMMA 1.   *Let $X_1, \ldots, X_n$ be i.i.d. real random variables with common density $f \in \mathbb{L}_\infty(\mathbb{R})$. Let $\mathbb{D}_1 = \mathbb{D}_3 = \mathbb{N}^*$ and $\mathbb{D}_2 = \{2^J, J \in \mathbb{N}\}$. For all $m = (l, D)$ with $l \in \{1, 2, 3\}$ and $D \in \mathbb{D}_l$, introduce $\{p_l, l \in \mathcal{L}_m\}$ defined as in (5.1) and*

$$Z_m = \frac{1}{n(n-1)}\sum_{i \neq j=1}^n H_m(X_i, X_j),$$

*with*

$$H_m(x, y) = \sum_{l \in \mathcal{L}_m}(p_l(x) - \langle f, p_l \rangle)(p_l(y) - \langle f, p_l \rangle).$$

*There exists some positive constant $C$ (depending only on $\varphi$) such that, for all $l \in \{1, 2, 3\}$, $D \in \mathbb{D}_l$, $x > 0$,*

$$(5.4) \quad \mathbb{P}\left(|Z_{(l,D)}| > \frac{C}{n}\left(\sqrt{Dx}(\|f\|_\infty + \sqrt{\|f\|_\infty}) + \|f\|_\infty x + \frac{Dx^2}{n}\right)\right) \leq 5.6e^{-x}.$$

The proof of this lemma is detailed in [9].

By setting $\lambda = \log(3/\beta)$ and $\tilde{\lambda} = \lambda + \log(5.6)$, Lemma 1 gives that there exists some constant $C > 0$ such that, for all $m = (l, D)$ in $\mathcal{M}$,

$$(5.5) \quad \mathbb{P}_f\left(U_n(H_m) < -\frac{C}{n}\left((\|f\|_\infty + \sqrt{\|f\|_\infty})\sqrt{D\tilde{\lambda}} + \|f\|_\infty\tilde{\lambda} + \frac{D\tilde{\lambda}^2}{n}\right)\right) \leq \frac{\beta}{3}.$$

We deduce from (5.3) and (5.5) that

$$
\begin{aligned}
&\mathbb{P}_f(T_\alpha \leq 0) \\
&\quad \leq \frac{\beta}{3} + \inf_{(l,D) \in \mathcal{M}} \Bigg\{ \mathbb{P}_f\Big( (P_n - P)(2\Pi_{S_{(l,D)}}(f) - 2f) + (P_n - P)(2f - 2f_0) \\
&\qquad\qquad + \|f - f_0\|_2^2 \leq \|f - \Pi_{S_{(l,D)}}(f)\|_2^2 + t_{(l,D)}(u_\alpha) \\
&\qquad\qquad + \frac{C}{n}\Big( (\|f\|_\infty + \sqrt{\|f\|_\infty})\sqrt{D\tilde{\lambda}} + \|f\|_\infty \tilde{\lambda} + \frac{D\tilde{\lambda}^2}{n} \Big) \Big) \Bigg\}.
\end{aligned}
$$

(5.6)

(b) Control of $(P_n - P)(2\Pi_{S_m}(f) - 2f)$ and $(P_n - P)(2f - 2f_0)$. We now use the following lemma due to Birgé and Massart [5], which provides a special version of Bernstein's inequality.

LEMMA 2. *Let $X_1, \ldots, X_n$ be independent random variables satisfying the moment condition*

$$
\frac{1}{n} \sum_{i=1}^n \mathbb{E}(|X_i|^k) \leq \frac{k!}{2} \nu c^{k-2} \qquad \text{for all } k \geq 2,
$$

*for some positive constants $\nu$ and $c$. Then, for any positive $x$,*

$$
\mathbb{P}\left( \frac{1}{n} \sum_{i=1}^n (X_i - \mathbb{E}(X_i)) \geq \frac{\sqrt{2\nu x}}{\sqrt{n}} + \frac{cx}{n} \right) \leq e^{-x}.
$$

*In particular, if for all $i$ in $\{1, \ldots, n\}$, $|X_i| \leq b$ and $\mathbb{E}(X_i^2) \leq \nu$, the above inequality is satisfied with $c = b/3$.*

It is easy to check that there exists some constant $C' > 0$ such that for all $l$ in $\{1, 2\}$, $D$ in $\mathbb{D}_l$,

$$
|2\Pi_{S_{(l,D)}}(f)(X_i) - 2f(X_i)| \leq C'\|f\|_\infty.
$$

Moreover, it is proved in [7], page 269, that one can take $C'$ such that for all $D$ in $\mathbb{D}_3$,

$$
|2\Pi_{S_{(3,D)}}(f)(X_i) - 2f(X_i)| \leq C'\|f\|_\infty \log(D+1).
$$

Since

$$
\mathbb{E}(2\Pi_{S_m}(f)(X_i) - 2f(X_i))^2 \leq 4\|f\|_\infty \|\Pi_{S_m}(f) - f\|_2^2,
$$

we can deduce from Lemma 2 that for all $m = (l, D) \in \mathcal{M}$,

$$
\begin{aligned}
\mathbb{P}_f\Big( (P_n - P)(2\Pi_{S_m}(f) - 2f) &< -2\sqrt{\|f\|_\infty} \|\Pi_{S_m}(f) - f\|_2 \sqrt{\frac{2\lambda}{n}} \\
&\quad - \frac{C'\lambda\|f\|_\infty \log(D+1)}{3n} \Big) \leq \frac{\beta}{3}.
\end{aligned}
$$

By using the elementary inequality $2ab \leq 4a^2/\varepsilon + \varepsilon b^2/4$, we obtain that for $m = (l, D) \in \mathcal{M}$,

$$
(5.7) \quad
\begin{aligned}
\mathbb{P}_f\Big( &(P_n - P)(2\Pi_{S_m}(f) - 2f) + \frac{\varepsilon}{4}\|\Pi_{S_m}(f) - f\|_2^2 \\
&< -\Big(\frac{8}{\varepsilon} + \frac{C'\log(D+1)}{3}\Big)\frac{\|f\|_\infty \lambda}{n}\Big) \leq \frac{\beta}{3}.
\end{aligned}
$$

The control of $(P_n - P)(2f - 2f_0)$ is computed in the same way and we get

$$
(5.8) \quad
\begin{aligned}
\mathbb{P}_f\Big( &(P_n - P)(2f - 2f_0) + \frac{\varepsilon}{4}\|f - f_0\|_2^2 \\
&< -\Big(2\Big(\frac{4}{\varepsilon} + \frac{1}{3}\Big)\|f\|_\infty + \frac{2}{3}\|f_0\|_\infty\Big)\frac{\lambda}{n}\Big) \leq \frac{\beta}{3}.
\end{aligned}
$$

Finally, we deduce from (5.6)–(5.8) that if there exists some $m = (l, D)$ in $\mathcal{M}$ such that

$$
\begin{aligned}
\Big(1 - \frac{\varepsilon}{4}\Big)\|f - f_0\|_2^2 >{} & \Big(1 + \frac{\varepsilon}{4}\Big)\|f - \Pi_{S_m}(f)\|_2^2 \\
&+ \frac{C}{n}\Big((\|f\|_\infty + \sqrt{\|f\|_\infty})\sqrt{D\tilde{\lambda}} + \|f\|_\infty \tilde{\lambda} + \frac{D\tilde{\lambda}^2}{n}\Big) \\
&+ \Big(\Big(\frac{16}{\varepsilon} + \frac{C'\log(D+1)+2}{3}\Big)\|f\|_\infty + \frac{2}{3}\|f_0\|_\infty\Big)\frac{\lambda}{n} \\
&+ t_m(u_\alpha),
\end{aligned}
$$

then

$$
\mathbb{P}_f(T_\alpha \leq 0) \leq \beta.
$$

This concludes the proof of Theorem 1.

5.2. *Proof of Corollary* 1. Assume that for all $l$ in $\{1, 2, 3\}$, $\mathcal{D}_l = \varnothing$ or $\{2^J, 0 \leq J \leq \log_2(n^2/(\log\log n)^3)\}$.

5.2.1. *An upper bound for* $t_m(u_\alpha)$, $m \in \mathcal{M}$.

PROPOSITION 1. *There exists some positive constant* $C(\alpha)$ *such that, for all* $m = (l, D)$ *in* $\mathcal{M}$,

$$
t_m(u_\alpha) \leq W_m(\alpha),
$$

*where*

$$
\begin{aligned}
W_m(\alpha) = \frac{C(\alpha)}{n}\Big( &(\|f_0\|_\infty + \sqrt{\|f_0\|_\infty})\sqrt{D\log\log n} + \frac{D(\log\log n)^2}{n} \\
&+ \|f_0\|_\infty(\log\log n)\log n\Big).
\end{aligned}
$$

PROOF. Recall that $t_m(u)$ denotes the $(1-u)$ quantile of the distribution of $\hat{T}_m$ under the hypothesis "$f = f_0$." We first notice that for $n \geq 16$, $|\mathcal{M}| \leq 3(1 + \log_2 n^2)$. So, setting $\alpha_n = \alpha/(3(1 + \log_2 n^2))$,

$$\mathbb{P}_{f_0}\bigg(\sup_{m \in \mathcal{M}} (\hat{T}_m - t_m(\alpha_n)) > 0\bigg) \leq \sum_{m \in \mathcal{M}} \mathbb{P}_{f_0}(\hat{T}_m - t_m(\alpha_n) > 0)$$

$$\leq \sum_{m \in \mathcal{M}} \frac{\alpha}{3(1 + \log_2 n^2)} \leq \alpha.$$

By definition of $u_\alpha$, this implies that $\alpha_n \leq u_\alpha$ and for all $m$ in $\mathcal{M}$,

$$t_m(u_\alpha) \leq t_m(\alpha_n).$$

It thus remains to give an upper bound for $t_m(\alpha_n)$. Let $m = (l, D) \in \mathcal{M}$. We use the same notation as in the proof of Theorem 1 to obtain the decomposition

$$\hat{T}_m = U_n(H_m) + (P_n - P)(2\Pi_{S_m}(f)) - 2P_n(f_0) + \|f_0\|_2^2 + \|\Pi_{S_m}(f)\|_2^2.$$

Under the null hypothesis "$f = f_0$,"

$$\hat{T}_m = U_n(H_m) + (P_n - P)(2\Pi_{S_m}(f_0) - 2f_0) - \|f_0\|_2^2 + \|\Pi_{S_m}(f_0)\|_2^2.$$

Since $\|\Pi_{S_m}(f_0) - f_0\|_2^2 = \|f_0\|_2^2 - \|\Pi_{S_m}(f_0)\|_2^2$, we obtain that, under "$f = f_0$,"

$$(5.9) \quad \hat{T}_m = U_n(H_m) + (P_n - P)(2\Pi_{S_m}(f_0) - 2f_0) - \|\Pi_{S_m}(f_0) - f_0\|_2^2.$$

As in the proof of Theorem 1, we control $U_n(H_m)$ from Lemma 1 and $(P_n - P)(2\Pi_{S_m}(f_0) - 2f_0)$ from Lemma 2. We set $\lambda_n = \log(2/\alpha_n)$ and $\tilde{\lambda}_n = \lambda_n + \log(5.6)$.

On one hand, Lemma 1 leads to

$$(5.10) \quad \mathbb{P}_{f_0}\bigg(U_n(H_m) > \frac{C}{n}\bigg(\sqrt{D\tilde{\lambda}_n}(\|f_0\|_\infty + \sqrt{\|f_0\|_\infty}) + \|f_0\|_\infty \tilde{\lambda}_n + \frac{D\tilde{\lambda}_n^2}{n}\bigg)\bigg) \leq \frac{\alpha_n}{2}.$$

On the other hand, since

$$(5.11) \quad |2(\Pi_{S_m}(f_0) - f_0)(X_i)| \leq C'\|f_0\|_\infty \log(D + 1)$$

and

$$\mathbb{E}_{f_0}(2(\Pi_{S_m}(f_0) - f_0)(X_i))^2 \leq 4\|f_0\|_\infty \|\Pi_{S_m}(f_0) - f_0\|_2^2,$$

it follows from Lemma 2 that

$$\mathbb{P}_{f_0}\bigg((P_n - P)(2\Pi_{S_m}(f_0) - 2f_0) > 2\sqrt{\|f_0\|_\infty}\|\Pi_{S_m}(f_0) - f_0\|_2\sqrt{\frac{2\lambda_n}{n}}$$
$$+ \frac{C'\|f_0\|_\infty\lambda_n\log(D+1)}{3n}\bigg) \leq \frac{\alpha_n}{2}.$$

Using the inequality $2ab \leq a^2 + b^2$, and the fact that for $n \geq 16$, $\log(D+1) \leq \log(n^2 + 1)$, we obtain that there exists $C'' > 0$ such that

$$(5.12) \quad \mathbb{P}_{f_0}\bigg((P_n - P)(2\Pi_{S_m}(f_0) - 2f_0) - \|\Pi_{S_m}(f_0) - f_0\|_2^2$$
$$> \frac{C''\|f_0\|_\infty\lambda_n\log n}{n}\bigg) \leq \frac{\alpha_n}{2}.$$

We derive from (5.9), (5.10) and (5.12) that

$$\mathbb{P}_{f_0}\bigg(\hat{T}_m > \frac{C}{n}\bigg(\sqrt{D\tilde{\lambda}_n}(\|f_0\|_\infty + \sqrt{\|f_0\|_\infty}) + \|f_0\|_\infty\tilde{\lambda}_n + \frac{D\tilde{\lambda}_n^2}{n}\bigg)$$
$$+ \frac{C''\|f_0\|_\infty\lambda_n\log n}{n}\bigg) \leq \alpha_n.$$

Finally, we notice that there exist some positive constants $c(\alpha)$ and $c'(\alpha)$ such that for $n \geq 3$, $\lambda_n \leq c(\alpha)\log\log n$ and $\tilde{\lambda}_n \leq c'(\alpha)\log\log n$, which completes the proof of Proposition 1. $\square$

5.2.2. *Uniform separation rates.* Let us fix $\beta$ in $]0,1[$ and $l$ in $\{1,2,3\}$ such that $\mathcal{D}_l = \{2^J, 0 \leq J \leq \log_2(n^2/(\log\log n)^3)\}$. From Theorem 1 and Proposition 1, we deduce that if $f$ satisfies

$$\|f - f_0\|_2^2 > (1 + \varepsilon)\inf_{D \in \mathcal{D}_l}\{\|f - \Pi_{S_{(l,D)}}(f)\|_2^2 + W_{(l,D)}(\alpha) + V_{(l,D)}(\beta)\},$$

then

$$\mathbb{P}_f(T_\alpha \leq 0) \leq \beta.$$

It is thus a matter of giving an upper bound for

$$\inf_{D \in \mathcal{D}_l}\{\|f - \Pi_{S_{(l,D)}}(f)\|_2^2 + W_{(l,D)}(\alpha) + V_{(l,D)}(\beta)\},$$

when $f$ belongs to some specified classes of functions. Recall that

$$\mathcal{B}_s^{(l)}(R, M)$$
$$= \{f \in \mathbb{L}_2(\mathbb{R}), \ \forall D \in \mathbb{D}_l, \ \|f - \Pi_{S_{(l,D)}}(f)\|_2^2 \leq R^2 D^{-2s}, \ \|f\|_\infty \leq M\}.$$

We now assume that $f$ belongs to $\mathcal{B}_s^{(l)}(R, M)$. Since $\|f - \Pi_{S_{(l,D)}}(f)\|_2^2 \le R^2 D^{-2s}$ and since the constant $C_2(\beta, \varepsilon, \|f\|_\infty, \|f_0\|_\infty)$ in Theorem 1 can be taken such that $C_2(\beta, \varepsilon, \|f\|_\infty, \|f_0\|_\infty) \le C_2(\beta, \varepsilon, M, \|f_0\|_\infty)$, we only need an upper bound for

$$\inf_{D \in \mathcal{D}_l} \left\{ R^2 D^{-2s} + C_1(\beta)(\sqrt{M} + M)\frac{\sqrt{D}}{n} \right.$$
$$+ C_1(\beta)\frac{D}{n^2} + C(\alpha)\frac{D(\log\log n)^2}{n^2}$$
$$+ C(\alpha)(\|f_0\|_\infty + \sqrt{\|f_0\|_\infty})\frac{\sqrt{D\log\log n}}{n}$$
$$\left. + C(\alpha)\|f_0\|_\infty \frac{(\log\log n)\log n}{n} + \frac{C_2(\beta, \varepsilon, M, \|f_0\|_\infty)}{n} \right\}.$$

Assuming that $n \ge 16$, this quantity is bounded from above by

$$\inf_{D \in \mathcal{D}_l} \left\{ R^2 D^{-2s} + (C_1(\beta) + C(\alpha))\frac{D(\log\log n)^2}{n^2} \right.$$
$$\left. + (C_1(\beta)(\sqrt{M} + M) + C(\alpha)(\|f_0\|_\infty + \sqrt{\|f_0\|_\infty}))\frac{\sqrt{D\log\log n}}{n} \right\}$$
$$+ (C(\alpha)\|f_0\|_\infty + C_2(\beta, \varepsilon, M, \|f_0\|_\infty))\frac{(\log\log n)\log n}{n}.$$

Since every $D$ in $\mathcal{D}_l$ is smaller than $n^2/(\log\log n)^3$,

$$\inf_{D \in \mathcal{D}_l} \left\{ R^2 D^{-2s} + \frac{\sqrt{D\log\log n}}{n} + \frac{D(\log\log n)^2}{n^2} \right\}$$
$$\le 2\inf_{D \in \mathcal{D}_l} \left\{ R^2 D^{-2s} + \frac{\sqrt{D\log\log n}}{n} \right\}.$$

We have $R^2 D^{-2s} < \sqrt{D\log\log n}/n$ if and only if $D > (R^4 n^2/\log\log n)^{1/(1+4s)}$, so we define $D_*$ by

$$\log_2(D_*) = [\log_2((R^4 n^2/(\log\log n))^{1/(1+4s)})] + 1,$$

and we consider three cases.

The first one is the case where $1 \le D_* \le 2^{[\log_2(n^2/(\log\log n)^3)]}$, which means that $D_* \in \mathcal{D}_l$. In this case, we have that

$$\inf_{D \in \mathcal{D}_l} \left\{ R^2 D^{-2s} + \frac{\sqrt{D\log\log n}}{n} \right\} \le R^2 D_*^{-2s} + \frac{\sqrt{D_*\log\log n}}{n}.$$

Since

$$R^2 D_*^{-2s} \le R^{2/(4s+1)} \left( \frac{\sqrt{\log\log n}}{n} \right)^{4s/(4s+1)}$$

and

$$\frac{\sqrt{D_* \log\log n}}{n} \leq \sqrt{2}\left(\frac{nR^2}{\sqrt{\log\log n}}\right)^{1/(4s+1)}\frac{\sqrt{\log\log n}}{n}$$

$$\leq \sqrt{2}R^{2/(4s+1)}\left(\frac{\sqrt{\log\log n}}{n}\right)^{4s/(4s+1)},$$

we obtain that

$$\inf_{D\in\mathcal{D}_l}\left\{R^2 D^{-2s} + \frac{\sqrt{D\log\log n}}{n}\right\} \leq (1+\sqrt{2})R^{2/(4s+1)}\left(\frac{\sqrt{\log\log n}}{n}\right)^{4s/(4s+1)}.$$

The second one is the case where $D_* > 2^{[\log_2(n^2/(\log\log n)^3)]}$. In this case, for all $D$ in $\mathcal{D}_l$,

$$\frac{\sqrt{D\log\log n}}{n} \leq R^2 D^{-2s}.$$

By taking $D_0 = 2^{[\log_2(n^2/(\log\log n)^3)]}$, we obtain that

$$\inf_{D\in\mathcal{D}_l}\left\{R^2 D^{-2s} + \frac{\sqrt{D\log\log n}}{n}\right\} \leq 2R^2 D_0^{-2s} \leq 2^{2s+1}R^2\left(\frac{(\log\log n)^3}{n^2}\right)^{2s}.$$

The third one is the case where $D_* < 1$. In this case, for all $D$ in $\mathcal{D}_l$, $R^2 D^{-2s} \leq \sqrt{D\log\log n}/n$, so by taking $D = 1$, we obtain that

$$\inf_{D\in\mathcal{D}_l}\left\{R^2 D^{-2s} + \frac{\sqrt{D\log\log n}}{n}\right\} \leq \frac{\sqrt{\log\log n}}{n}.$$

This completes the proof of the corollary.

5.3. *Proof of Corollary* 2. We use the same notation as in Theorem 1. We assume that $f_0 = \mathbb{I}_{[0,1]}$ and $\mathcal{M} = \{(1, D), D \in \mathcal{D}_1\}$, which means that we only consider spaces generated by constant piecewise functions. We first prove that if there exists $m$ in $\mathcal{M}$ such that

(5.13)     $$\|\Pi_{S_m}(f) - f_0\|_2^2 \geq (1+\varepsilon)(t_m(u_\alpha) + V_m(\beta)),$$

then

$$\mathbb{P}_f(T_\alpha \leq 0) \leq \beta.$$

Using inequality (5.2) and the definition of $\hat{m}$, we derive that

$$\mathbb{P}_f(T_\alpha \leq 0) \leq \inf_{m\in\mathcal{M}} \mathbb{P}_f\bigg(U_n(H_m) + (P_n - P)(2\Pi_{S_m}(f) - 2f_0)$$

$$+ \|\Pi_{S_m}(f) - f_0\|_2^2 + 2\int_0^1 f_0(f - \Pi_{S_m}(f)) \leq t_m(u_\alpha)\bigg).$$

Since $f_0 = \mathbb{I}_{[0,1]}$ and since, for all $m$ in $\mathcal{M}$, $S_m$ is generated by constant piecewise functions, we have

$$\int_0^1 f_0(f - \Pi_{S_m}(f)) = 0.$$

We then obtain that

$$\mathbb{P}_f(T_\alpha \leq 0) \leq \inf_{m \in \mathcal{M}} \mathbb{P}_f(U_n(H_m) + (P_n - P)(2\Pi_{S_m}(f) - 2f_0)$$
$$+ \|\Pi_{S_m}(f) - f_0\|_2^2 \leq t_m(u_\alpha)),$$

and we conclude as in the proof of Theorem 1.

Let for all $s > 0$ and $R > 0$, $\mathcal{H}_s(R)$ be defined by (2.4). Ingster [14] proved that for all functions $g$ in $\mathcal{H}_s(R)$, for all $D$ in $\mathbb{N}^*$,

$$\|\Pi_{S_{(1,D)}}(g)\|_2^2 \geq C_1 \|g\|_2^2 - C_2 D^{-2s},$$

where $C_1 \in ]0, 1[$ and $C_2 > 0$ (see [14], part III, inequality (5.16)).

Assume that $f \in \mathcal{H}_s(R)$. Since $f_0 = \mathbb{I}_{[0,1]}$, $f - f_0 \in \mathcal{H}_s(R)$. We then derive from Ingster's inequality that for all $D$ in $\mathcal{D}_1$,

$$\|\Pi_{S_{(1,D)}}(f - f_0)\|_2^2 \geq C_1 \|f - f_0\|_2^2 - C_2 D^{-2s}.$$

Moreover, $\Pi_{S_{(1,D)}}(f - f_0) = \Pi_{S_{(1,D)}}(f) - f_0$. As a consequence, if there exists $m = (1, D)$ in $\mathcal{M}$ such that

$$C_1 \|f - f_0\|_2^2 \geq C_2 D^{-2s} + (1 + \varepsilon)(t_m(u_\alpha) + V_m(\beta)),$$

then (5.13) holds. The end of the proof is similar to the proof of Corollary 1.

5.4. *Proof of Theorem* 2. Let $\tilde{T}_\alpha$ be the test statistic defined by (3.4). We have that

$$\mathbb{P}_f(\tilde{T}_\alpha \leq 0) \leq \inf_{m \in \mathcal{M}} \mathbb{P}_f(\tilde{T}_m(X_1, \ldots, X_n) \leq \tilde{q}_{m,\alpha}),$$

and since $\mathcal{D}_2 \neq \varnothing$, setting $\mathcal{M}' = \{(2, D), \ D \in \mathcal{D}_2\}$,

(5.14) $$\mathbb{P}_f(\tilde{T}_\alpha \leq 0) \leq \inf_{m \in \mathcal{M}'} \mathbb{P}_f(\tilde{T}_m(X_1, \ldots, X_n) \leq \tilde{q}_{m,\alpha}).$$

Let us introduce some notation. Using the definitions given in (5.1), for $m$ in $\mathcal{M}'$, we recall that

$$\hat{T}_m\left(\frac{X_1 - \mu}{\sigma}, \ldots, \frac{X_n - \mu}{\sigma}\right) = \frac{1}{n(n-1)} \sum_{l \in \mathcal{L}_m} \sum_{i \neq j=1}^n p_l\left(\frac{X_i - \mu}{\sigma}\right) p_l\left(\frac{X_j - \mu}{\sigma}\right)$$
$$+ \|f_0\|_2^2 - \frac{2}{n} \sum_{i=1}^n f_0\left(\frac{X_i - \mu}{\sigma}\right),$$

and we set

$$a_l(\mu, \sigma) = \int p_l\left(\frac{x - \mu}{\sigma}\right) f(x)\, dx,$$

$$b(\mu, \sigma) = \int f_0\left(\frac{x - \mu}{\sigma}\right) f(x)\, dx,$$

$$Z_m(\mu, \sigma) = \frac{1}{n(n-1)} \sum_{l \in \mathcal{L}_m} \sum_{i \neq j = 1}^{n} \left[ p_l\left(\frac{X_i - \mu}{\sigma}\right) - a_l(\mu, \sigma) \right]$$
$$\times \left[ p_l\left(\frac{X_j - \mu}{\sigma}\right) - a_l(\mu, \sigma) \right],$$

$$L_m(\mu, \sigma) = (P_n - P)\left( 2 \sum_{l \in \mathcal{L}_m} a_l(\mu, \sigma) p_l\left(\frac{\cdot - \mu}{\sigma}\right) - 2 f_0\left(\frac{\cdot - \mu}{\sigma}\right) \right).$$

The following identity holds:

$$\hat{T}_m\left(\frac{X_1 - \mu}{\sigma}, \ldots, \frac{X_n - \mu}{\sigma}\right) = Z_m(\mu, \sigma) + L_m(\mu, \sigma) + \|f_0\|_2^2$$
$$+ \sum_{l \in \mathcal{L}_m} a_l^2(\mu, \sigma) - 2b(\mu, \sigma).$$

Recall that $S_m$ is the linear subspace of $\mathbb{L}_2(\mathbb{R})$ spanned by the functions $\{p_l, l \in \mathcal{L}_m\}$ and that $\Pi_{S_m}$ denotes the orthogonal projection onto $S_m$ in $\mathbb{L}_2(\mathbb{R})$. Then

$$(5.15) \qquad \sum_{l \in \mathcal{L}_m} a_l^2(\mu, \sigma) = \|\Pi_{S_m}(\sigma f(\sigma \cdot + \mu))\|_2^2.$$

After some easy computations, one can prove that

$$\|f_0\|_2^2 + \sum_{l \in \mathcal{L}_m} a_l^2(\mu, \sigma) - 2b(\mu, \sigma)$$
$$= \sigma \left\| \frac{1}{\sigma} f_0\left(\frac{\cdot - \mu}{\sigma}\right) - f \right\|_2^2$$
$$- \|\sigma f(\sigma \cdot + \mu) - \Pi_{S_m}(\sigma f(\sigma \cdot + \mu))\|_2^2.$$

This implies that

$$\mathbb{P}_f(\tilde{T}_m(X_1, \ldots, X_n) \leq \tilde{q}_{m,\alpha})$$
$$= \mathbb{P}_f\left( \inf_{(\mu, \sigma) \in K} \left\{ Z_m(\mu, \sigma) + L_m(\mu, \sigma) \right. \right.$$
$$+ \sigma \left\| \frac{1}{\sigma} f_0\left(\frac{\cdot - \mu}{\sigma}\right) - f \right\|_2^2$$
$$\left. \left. - \|\sigma f(\sigma \cdot + \mu) - \Pi_{S_m}(\sigma f(\sigma \cdot + \mu))\|_2^2 \right\} \leq \tilde{q}_{m,\alpha} \right),$$

or

$$\mathbb{P}_f(\tilde{T}_m(X_1,\ldots,X_n) \le \tilde{q}_{m,\alpha})$$

$$= \mathbb{P}_f\Bigg(\inf_{(\mu,\sigma)\in K}\Bigg\{Z_m(\mu,\sigma) + L_m(\mu,\sigma) + \frac{\varepsilon}{4}\sigma\left\|\frac{1}{\sigma}f_0\left(\frac{\cdot-\mu}{\sigma}\right)-f\right\|_2^2$$

$$+ \frac{\varepsilon}{4}\|\sigma f(\sigma\cdot+\mu) - \Pi_{S_m}(\sigma f(\sigma\cdot+\mu))\|_2^2$$

$$+ \left(1-\frac{\varepsilon}{4}\right)\sigma\left\|\frac{1}{\sigma}f_0\left(\frac{\cdot-\mu}{\sigma}\right)-f\right\|_2^2$$

$$- \left(1+\frac{\varepsilon}{4}\right)\|\sigma f(\sigma\cdot+\mu) - \Pi_{S_m}(\sigma f(\sigma\cdot+\mu))\|_2^2\Bigg\} \le \tilde{q}_{m,\alpha}\Bigg).$$

Therefore,

$$\mathbb{P}_f(\tilde{T}_m \le \tilde{q}_{m,\alpha})$$

$$\le \mathbb{P}_f\Bigg(\sup_{(\mu,\sigma)\in K}\Bigg\{-Z_m(\mu,\sigma) - L_m(\mu,\sigma)$$

$$- \frac{\varepsilon}{4}\sigma\left\|\frac{1}{\sigma}f_0\left(\frac{\cdot-\mu}{\sigma}\right)-f\right\|_2^2$$

$$- \frac{\varepsilon}{4}\|\sigma f(\sigma\cdot+\mu) - \Pi_{S_m}(\sigma f(\sigma\cdot+\mu))\|_2^2\Bigg\}$$

$$\ge \left(1-\frac{\varepsilon}{4}\right)\inf_{(\mu,\sigma)\in K}\sigma\left\|\frac{1}{\sigma}f_0\left(\frac{\cdot-\mu}{\sigma}\right)-f\right\|_2^2$$

$$- \left(1+\frac{\varepsilon}{4}\right)\sup_{(\mu,\sigma)\in K}\|\sigma f(\sigma\cdot+\mu) - \Pi_{S_m}(\sigma f(\sigma\cdot+\mu))\|_2^2 - \tilde{q}_{m,\alpha}\Bigg).$$

Let $\tau_m(\beta)$ denote the $(1-\beta)$ quantile of $\sup_{(\mu,\sigma)\in K}\Gamma_m(\mu,\sigma)$, where

$$\Gamma_m(\mu,\sigma) = -Z_m(\mu,\sigma) - L_m(\mu,\sigma) - \frac{\varepsilon}{4}\sigma\left\|\frac{1}{\sigma}f_0\left(\frac{\cdot-\mu}{\sigma}\right)-f\right\|_2^2$$

$$- \frac{\varepsilon}{4}\|\sigma f(\sigma\cdot+\mu) - \Pi_{S_m}(\sigma f(\sigma\cdot+\mu))\|_2^2.$$

It follows from the above inequality and from (5.14) that if

$$\left(1-\frac{\varepsilon}{4}\right)\inf_{(\mu,\sigma)\in K}\sigma\left\|\frac{1}{\sigma}f_0\left(\frac{\cdot-\mu}{\sigma}\right)-f\right\|_2^2$$

$$> \inf_{m\in\mathcal{M}'}\Bigg(\left(1+\frac{\varepsilon}{4}\right)\sup_{(\mu,\sigma)\in K}\|\sigma f(\sigma\cdot+\mu) - \Pi_{S_m}(\sigma f(\sigma\cdot+\mu))\|_2^2$$

$$+ \tau_m(\beta) + \tilde{q}_{m,\alpha}\Bigg),$$

then

$$\mathbb{P}_f(\tilde{T}_\alpha \leq 0) \leq \beta.$$

Let $m = (2, 2^J) \in \mathcal{M}'$. We now need to compute an upper bound for $\tau_m(\beta)$.

To compute this upper bound, we introduce a finite grid on $K = [\underline{\mu}, \overline{\mu}] \times [\underline{\sigma}, \overline{\sigma}]$.

Let $\mu_0 < \mu_1 < \cdots < \mu_N$ and $\sigma_0 < \sigma_1 < \cdots < \sigma_N$ such that $\mu_0 = \underline{\mu}$, $\mu_N = \overline{\mu}$, $\sigma_0 = \underline{\sigma}$, $\sigma_N = \overline{\sigma}$, for all $\delta$ in $\{0, \ldots, N-2\}$, $|\mu_{\delta+1} - \mu_\delta| = \Delta_\mu$, $|\sigma_{\delta+1} - \sigma_\delta| = \Delta_\sigma$ and $|\mu_N - \mu_{N-1}| \leq \Delta_\mu$, $|\sigma_N - \sigma_{N-1}| \leq \Delta_\sigma$. $\Delta_\mu$ and $\Delta_\sigma$ will be chosen later. Let for $(\delta, \delta')$ in $\{0, \ldots, N-1\}^2$,

$$A_{\delta, \delta'} = [\mu_\delta, \mu_{\delta+1}] \times [\sigma_{\delta'}, \sigma_{\delta'+1}].$$

The following inequality holds:

$$\sup_{(\mu, \sigma) \in K} \Gamma_m(\mu, \sigma) \leq \sup_{(\delta, \delta') \in \{0, \ldots, N-1\}^2} \Gamma_m(\mu_\delta, \sigma_{\delta'})$$
$$+ \sup_{(\delta, \delta') \in \{0, \ldots, N-1\}^2} \sup_{(\mu, \sigma) \in A_{\delta, \delta'}} (\Gamma_m(\mu, \sigma) - \Gamma_m(\mu_\delta, \sigma_{\delta'})).$$

5.4.1. *Control of* $\sup_{(\delta, \delta') \in \{0, \ldots, N-1\}^2} \Gamma_m(\mu_\delta, \sigma_{\delta'})$. Introduce

$$L_m^{(1)}(\mu, \sigma) = (P_n - P)\left(2\sigma f - 2(\Pi_{S_m}(\sigma f(\sigma. + \mu)))\left(\frac{\cdot - \mu}{\sigma}\right)\right)$$
$$- \frac{\varepsilon}{4}\|\sigma f(\sigma. + \mu) - \Pi_{S_m}(\sigma f(\sigma. + \mu))\|_2^2$$

and

$$L_m^{(2)}(\mu, \sigma) = (P_n - P)\left(2f_0\left(\frac{\cdot - \mu}{\sigma}\right) - 2\sigma f\right) - \frac{\varepsilon}{4}\sigma\left\|\frac{1}{\sigma}f_0\left(\frac{\cdot - \mu}{\sigma}\right) - f\right\|_2^2.$$

Since

$$-L_m(\mu, \sigma) = (P_n - P)\left(2\sigma f - 2(\Pi_{S_m}(\sigma f(\sigma. + \mu)))\left(\frac{\cdot - \mu}{\sigma}\right)\right)$$
$$+ (P_n - P)\left(2f_0\left(\frac{\cdot - \mu}{\sigma}\right) - 2\sigma f\right),$$

we have that

$$\sup_{(\delta, \delta') \in \{0, \ldots, N-1\}^2} \Gamma_m(\mu_\delta, \sigma_{\delta'})$$
$$\leq \sup_{(\delta, \delta') \in \{0, \ldots, N-1\}^2} |Z_m(\mu_\delta, \sigma_{\delta'})|$$
$$+ \sup_{(\delta, \delta') \in \{0, \ldots, N-1\}^2} L_m^{(1)}(\mu_\delta, \sigma_{\delta'}) + \sup_{(\delta, \delta') \in \{0, \ldots, N-1\}^2} L_m^{(2)}(\mu_\delta, \sigma_{\delta'}).$$

(a) Control of $\sup_{(\delta,\delta')\in\{0,\ldots,N-1\}^2}|Z_m(\mu_\delta,\sigma_{\delta'})|$. We apply Lemma 1 by replacing the $X_i$'s by the variables $(X_i-\mu)/\sigma$ and the density $f$ by the density of the $(X_i-\mu)/\sigma$'s that is $\sigma f(\sigma.+\mu)$. With $m=(2,2^J)$ and $\{p_l, l\in \mathcal{L}_m\}=\{\varphi_{J,k}, k\in\mathbb{Z}\}$, we obtain that there exists some constant $C>0$ such that for any $(\mu,\sigma)$ in $K$, for any $x>0$,

$$\mathbb{P}_f\left(|Z_m(\mu,\sigma)|>\frac{C}{n}\left(2^{J/2}\sqrt{x}(\overline{\sigma}\|f\|_\infty+\sqrt{\overline{\sigma}\|f\|_\infty})+\overline{\sigma}\|f\|_\infty x+\frac{2^J}{n}x^2\right)\right)$$
$$\leq 5.6e^{-x}.$$

Hence,

$$\begin{aligned}
\mathbb{P}_f&\left(\sup_{(\delta,\delta')\in\{0,\ldots,N-1\}^2}|Z_m(\mu_\delta,\sigma_{\delta'})|\right.\\
(5.16)\qquad&>\frac{C}{n}\left(2^{J/2}\sqrt{x}(\overline{\sigma}\|f\|_\infty+\sqrt{\overline{\sigma}\|f\|_\infty})+\overline{\sigma}\|f\|_\infty x+\frac{2^J}{n}x^2\right)\right)\\
&\leq 5.6N^2e^{-x}.
\end{aligned}$$

(b) Control of $\sup_{(\delta,\delta')\in\{0,\ldots,N-1\}^2}L_m^{(1)}(\mu_\delta,\sigma_{\delta'})$. To get an upper bound for this supremum, we use Lemma 2. We recall that for $m=(2,2^J)$, $\{p_l, l\in \mathcal{L}_m\}=\{\varphi_{J,k}, k\in\mathbb{Z}\}$, where $\varphi_{J,k}=2^{J/2}\varphi(2^J\cdot-k)$ and $\varphi$ is a compactly supported scaling function. Then, for all $l$ in $\mathcal{L}_m$,

$$\|p_l\|_\infty\leq 2^{J/2}\|\varphi\|_\infty.$$

Assume that the support of the function $\varphi$ is included in $[-\Phi,\Phi]$, with $\Phi>0$. This implies that for any $x$ in $\mathbb{R}$, the cardinality of the set $\{l\in\mathcal{L}_m, p_l(x)\neq 0\}$ is smaller than $2\Phi$ and that for all $(\mu,\sigma)$ in $K$, $l$ in $\mathcal{L}_m$,

$$|a_l(\mu,\sigma)|\leq 2\Phi\sigma 2^{-J/2}\|f\|_\infty\|\varphi\|_\infty.$$

Hence, for all $x$ in $\mathbb{R}$,

$$\begin{aligned}
\left|(\Pi_{S_m}(\sigma f(\sigma.+\mu)))\left(\frac{x-\mu}{\sigma}\right)\right|&=\left|\sum_{l\in\mathcal{L}_m}a_l(\mu,\sigma)p_l\left(\frac{x-\mu}{\sigma}\right)\right|\\
&\leq 2\Phi\sigma 2^{-J/2}\|f\|_\infty\|\varphi\|_\infty\sum_{l\in\mathcal{L}_m}\left|p_l\left(\frac{x-\mu}{\sigma}\right)\right|\\
&\leq 4\Phi^2\overline{\sigma}\|f\|_\infty\|\varphi\|_\infty^2.
\end{aligned}$$

Therefore, there exists some constant $C'>0$ such that

$$\left\|2\sigma f-2(\Pi_{S_m}(\sigma f(\sigma.+\mu)))\left(\frac{\cdot-\mu}{\sigma}\right)\right\|_\infty\leq C'\overline{\sigma}\|f\|_\infty.$$

Furthermore, we have

$$\int \left(2\sigma f - 2(\Pi_{S_m}(\sigma f(\sigma . + \mu)))\left(\frac{\cdot - \mu}{\sigma}\right)\right)^2 f$$
$$\leq 4\overline{\sigma}\|f\|_\infty \|\sigma f(\sigma . + \mu) - \Pi_{S_m}(\sigma f(\sigma . + \mu))\|_2^2.$$

Then, for all $(\mu, \sigma)$ in $K$, Lemma 2 gives the following exponential inequality: for all $x > 0$,

$$\mathbb{P}_f\Big((P_n - P)\Big(2\sigma f - 2(\Pi_{S_m}(\sigma f(\sigma . + \mu)))\left(\frac{\cdot - \mu}{\sigma}\right)\Big)$$
$$> \frac{C'\overline{\sigma}\|f\|_\infty}{3}\frac{x}{n}$$
$$+ 2\sqrt{2\overline{\sigma}\|f\|_\infty}\|\sigma f(\sigma . + \mu) - \Pi_{S_m}(\sigma f(\sigma . + \mu))\|_2\sqrt{\frac{x}{n}}\Big) \leq e^{-x}.$$

By using the elementary inequality $2ab \leq 4a^2/\varepsilon + \varepsilon b^2/4$, we obtain that

$$(5.17)\ \mathbb{P}_f\Big(\sup_{(\delta, \delta') \in \{0, \dots, N-1\}^2} L_m^{(1)}(\mu_\delta, \sigma_{\delta'}) > \Big(\frac{8}{\varepsilon} + \frac{C'}{3}\Big)\overline{\sigma}\|f\|_\infty \frac{x}{n}\Big) \leq N^2 e^{-x}.$$

(c) Control of $\sup_{(\delta, \delta') \in \{0, \dots, N-1\}^2} L_m^{(2)}(\mu_\delta, \sigma_{\delta'})$. We control $\sup_{(\delta, \delta') \in \{0, \dots, N-1\}^2} L_m^{(2)}(\mu_\delta, \sigma_{\delta'})$ in the same way, noticing that for any $(\mu, \sigma)$ in $K$,

$$\left\|2f_0\left(\frac{\cdot - \mu}{\sigma}\right) - 2\sigma f\right\|_\infty \leq 2(\|f_0\|_\infty + \overline{\sigma}\|f\|_\infty)$$

and

$$\int \left(2f_0\left(\frac{\cdot - \mu}{\sigma}\right) - 2\sigma f\right)^2 f \leq 4\sigma^2 \|f\|_\infty \left\|\frac{1}{\sigma}f_0\left(\frac{\cdot - \mu}{\sigma}\right) - f\right\|_2^2.$$

We get

$$(5.18)\ \begin{aligned}&\mathbb{P}_f\Big(\sup_{(\delta, \delta') \in \{0, \dots, N-1\}^2} L_m^{(2)}(\mu_\delta, \sigma_{\delta'})\\ &> 2\Big(\frac{4}{\varepsilon}\overline{\sigma}\|f\|_\infty + \frac{\|f_0\|_\infty + \overline{\sigma}\|f\|_\infty}{3}\Big)\frac{x}{n}\Big) \leq N^2 e^{-x}.\end{aligned}$$

Assume that $N \geq 2$. Collecting the inequalities (5.16), (5.17) and (5.18) and choosing $x$ such that $7.6N^2 e^{-x} = \beta/2$, we obtain that there exists some positive constant $C = C(\overline{\sigma}, \varepsilon, \|f\|_\infty, \|f_0\|_\infty, \beta)$ such that

$$(5.19)\ \begin{aligned}&\mathbb{P}_f\Big(\sup_{(\delta, \delta') \in \{0, \dots, N-1\}^2} \Gamma_m(\mu_\delta, \sigma_{\delta'})\\ &> C\Big(2^{J/2}\frac{\sqrt{\log N}}{n} + \frac{\log N}{n} + 2^J\frac{\log^2 N}{n^2}\Big)\Big) \leq \frac{\beta}{2}.\end{aligned}$$

5.4.2. *Control of* $\sup_{(\delta,\delta')\in\{0,\ldots,N-1\}^2,(\mu,\sigma)\in A_{\delta,\delta'}}(\Gamma_m(\mu,\sigma)-\Gamma_m(\mu_\delta,\sigma_{\delta'}))$. We have

$$\Gamma_m(\mu,\sigma) = -\frac{1}{n(n-1)}\sum_{l\in\mathcal{L}_m}\sum_{i\neq j=1}^n p_l\left(\frac{X_i-\mu}{\sigma}\right)p_l\left(\frac{X_j-\mu}{\sigma}\right)$$

$$+ \sum_{l\in\mathcal{L}_m}a_l^2(\mu,\sigma) - \frac{\varepsilon}{4}\|\sigma f(\sigma\cdot+\mu)-\Pi_{S_m}(\sigma f(\sigma\cdot+\mu))\|_2^2$$

$$- \frac{\varepsilon}{4}\sigma\left\|\frac{1}{\sigma}f_0\left(\frac{\cdot-\mu}{\sigma}\right)-f\right\|_2^2 + (P_n-P)\left(2f_0\left(\frac{\cdot-\mu}{\sigma}\right)\right).$$

Using (5.15), this implies that

$$\Gamma_m(\mu,\sigma) = -\frac{1}{n(n-1)}\sum_{l\in\mathcal{L}_m}\sum_{i\neq j=1}^n p_l\left(\frac{X_i-\mu}{\sigma}\right)p_l\left(\frac{X_j-\mu}{\sigma}\right)$$

$$- \frac{\varepsilon}{2}\sigma\|f\|_2^2 - \frac{\varepsilon}{4}\|f_0\|_2^2 + \left(1+\frac{\varepsilon}{4}\right)\sum_{l\in\mathcal{L}_m}a_l^2(\mu,\sigma) + \frac{\varepsilon}{2}b(\mu,\sigma)$$

$$+ (P_n-P)\left(2f_0\left(\frac{\cdot-\mu}{\sigma}\right)\right).$$

Let

$$\Gamma_m^{(1)}(\mu,\sigma,\delta,\delta') = -\frac{1}{n(n-1)}$$

$$\times \sum_{l\in\mathcal{L}_m}\sum_{i\neq j=1}^n\left(p_l\left(\frac{X_i-\mu}{\sigma}\right)p_l\left(\frac{X_j-\mu}{\sigma}\right)\right.$$

$$\left. - p_l\left(\frac{X_i-\mu_\delta}{\sigma_{\delta'}}\right)p_l\left(\frac{X_j-\mu_\delta}{\sigma_{\delta'}}\right)\right),$$

$$\Gamma_m^{(2)}(\mu,\sigma,\delta,\delta') = \sum_{l\in\mathcal{L}_m}a_l^2(\mu,\sigma) - \sum_{l\in\mathcal{L}_m}a_l^2(\mu_\delta,\sigma_{\delta'}),$$

$$\Gamma_m^{(3)}(\mu,\sigma,\delta,\delta') = b(\mu,\sigma) - b(\mu_\delta,\sigma_{\delta'}),$$

$$\Gamma_m^{(4)}(\mu,\sigma,\delta,\delta') = (P_n-P)\left(2f_0\left(\frac{\cdot-\mu}{\sigma}\right) - 2f_0\left(\frac{\cdot-\mu_\delta}{\sigma_{\delta'}}\right)\right).$$

Then

$$\sup_{(\delta,\delta')\in\{0,\ldots,N-1\}^2}\sup_{(\mu,\sigma)\in A_{\delta,\delta'}}(\Gamma_m(\mu,\sigma)-\Gamma_m(\mu_\delta,\sigma_{\delta'}))$$

$$\leq \sup_{(\delta,\delta')}\sup_{(\mu,\sigma)\in A_{\delta,\delta'}}|\Gamma_m^{(1)}(\mu,\sigma,\delta,\delta')|$$

$$+ \left(1 + \frac{\varepsilon}{4}\right) \sup_{(\delta,\delta')} \sup_{(\mu,\sigma)\in A_{\delta,\delta'}} |\Gamma_m^{(2)}(\mu,\sigma,\delta,\delta')|$$

$$+ \frac{\varepsilon}{2} \sup_{(\delta,\delta')} \sup_{(\mu,\sigma)\in A_{\delta,\delta'}} |\Gamma_m^{(3)}(\mu,\sigma,\delta,\delta')| + \sup_{(\delta,\delta')} \sup_{(\mu,\sigma)\in A_{\delta,\delta'}} |\Gamma_m^{(4)}(\mu,\sigma,\delta,\delta')|$$

$$+ \frac{\varepsilon}{2} \Delta_\sigma \|f\|_\infty.$$

(a) Control of $\sup_{(\delta,\delta')\in\{0,\ldots,N-1\}^2} \sup_{(\mu,\sigma)\in A_{\delta,\delta'}} |\Gamma_m^{(1)}(\mu,\sigma,\delta,\delta')|.$   One can easily see that

$$\Gamma_m^{(1)}(\mu,\sigma,\delta,\delta')$$
$$= -\frac{1}{n(n-1)}$$
$$\times \sum_{l\in\mathcal{L}_m} \sum_{i\neq j=1}^n \left[ \left(p_l\left(\frac{X_i-\mu}{\sigma}\right) - p_l\left(\frac{X_i-\mu_\delta}{\sigma_{\delta'}}\right)\right) p_l\left(\frac{X_j-\mu}{\sigma}\right) \right.$$
$$\left. + \left(p_l\left(\frac{X_j-\mu}{\sigma}\right) - p_l\left(\frac{X_j-\mu_\delta}{\sigma_{\delta'}}\right)\right) p_l\left(\frac{X_i-\mu_\delta}{\sigma_{\delta'}}\right) \right].$$

We recall that for all $x,y$ in $\mathbb{R}$, the cardinality of the set $\{l\in\mathcal{L}_m, p_l(x)\neq p_l(y)\}$ is not larger than $4\Phi$ and that for $l$ in $\mathcal{L}_m$,

$$\|p_l\|_\infty \leq 2^{J/2}\|\varphi\|_\infty.$$

Since $\varphi$ is a Lipschitz function with Lipschitz constant $C_\varphi$, we also have that for all $x,y$ in $\mathbb{R}$,

$$(5.20) \qquad \sum_{l\in\mathcal{L}_m} |p_l(x) - p_l(y)| \leq 4\Phi 2^{3J/2} C_\varphi |x-y|.$$

This implies that

$$|\Gamma_m^{(1)}(\mu,\sigma,\delta,\delta')| \leq 8\Phi C_\varphi \|\varphi\|_\infty \frac{2^{2J}}{n} \sum_{i=1}^n \left|\frac{X_i-\mu}{\sigma} - \frac{X_i-\mu_\delta}{\sigma_{\delta'}}\right|$$

and

$$\sup_{(\mu,\sigma)\in A_{\delta,\delta'}} |\Gamma_m^{(1)}(\mu,\sigma,\delta,\delta')|$$
$$\leq 8\Phi C_\varphi \|\varphi\|_\infty 2^{2J} \left\{ \frac{1}{n}\sum_{i=1}^n |X_i|\left(\frac{\Delta_\sigma}{\underline{\sigma}^2}\right) + \frac{1}{\underline{\sigma}^2}(\overline{\sigma}\Delta_\mu + (|\overline{\mu}|\vee|\underline{\mu}|)\Delta_\sigma) \right\}.$$

Assuming that $(h_2)$ holds, we derive from Lemma 2 that for all $x > 0$,

$$(5.21) \qquad \mathbb{P}_f\left(\sum_{i=1}^n (|X_i| - \mathbb{E}(|X_i|)) \geq \sqrt{2\nu n x} + cx\right) \leq e^{-x}.$$

Since $\mathbb{E}(|X_i|) \leq \sqrt{\nu}$ and $N \geq 2$, by taking $x = \log(4N^2/\beta)$, we obtain that there exists some constant $C(\underline{\mu}, \overline{\mu}, \underline{\sigma}, \overline{\sigma}, \nu, c, \beta) > 0$ such that

$$
(5.22) \quad
\begin{aligned}
\mathbb{P}_f \Big( \sup_{(\delta, \delta') \in \{0, \ldots, N-1\}^2} &\sup_{(\mu, \sigma) \in A_{\delta, \delta'}} |\Gamma_m^{(1)}(\mu, \sigma, \delta, \delta')| \\
&\geq C(\underline{\mu}, \overline{\mu}, \underline{\sigma}, \overline{\sigma}, \nu, c, \beta) 2^{2J} (\Delta_\mu \vee \Delta_\sigma) \log N \Big) \leq \frac{\beta}{4}.
\end{aligned}
$$

(b) Control of $\sup_{(\delta, \delta') \in \{0, \ldots, N-1\}^2} \sup_{(\mu, \sigma) \in A_{\delta, \delta'}} |\Gamma_m^{(2)}(\mu, \sigma, \delta, \delta')|$. Since

$$
\Gamma_m^{(2)}(\mu, \sigma, \delta, \delta') = \sum_{l \in \mathcal{L}_m} (a_l(\mu, \sigma) - a_l(\mu_\delta, \sigma_{\delta'}))(a_l(\mu, \sigma) + a_l(\mu_\delta, \sigma_{\delta'}))
$$

and $|a_l(\mu, \sigma) + a_l(\mu_\delta, \sigma_{\delta'})| \leq 2.2^{J/2} \|\varphi\|_\infty$,

$$
|\Gamma_m^{(2)}(\mu, \sigma, \delta, \delta')| \leq 2.2^{J/2} \|\varphi\|_\infty \int \sum_{l \in \mathcal{L}_m} \left| p_l \left( \frac{x - \mu}{\sigma} \right) - p_l \left( \frac{x - \mu_\delta}{\sigma_{\delta'}} \right) \right| f(x) \, dx.
$$

By (5.20), we obtain

$$
\begin{aligned}
&|\Gamma_m^{(2)}(\mu, \sigma, \delta, \delta')| \\
&\leq 8\Phi \|\varphi\|_\infty C_\varphi 2^{2J} \int \left| \frac{x - \mu}{\sigma} - \frac{x - \mu_\delta}{\sigma_{\delta'}} \right| f(x) \, dx, \\
&\leq 8\Phi \|\varphi\|_\infty C_\varphi 2^{2J} \int \left( |x| \frac{\Delta_\sigma}{\underline{\sigma}^2} + \frac{1}{\underline{\sigma}^2} (\overline{\sigma} \Delta_\mu + (|\overline{\mu}| \vee |\underline{\mu}|) \Delta_\sigma) \right) f(x) \, dx.
\end{aligned}
$$

We deduce from $(h_2)$ that there exists some constant $C(\underline{\mu}, \overline{\mu}, \underline{\sigma}, \overline{\sigma}, \nu) > 0$ such that

$$
(5.23) \quad
\begin{aligned}
\sup_{(\delta, \delta') \in \{0, \ldots, N-1\}^2} &\sup_{(\mu, \sigma) \in A_{\delta, \delta'}} |\Gamma_m^{(2)}(\mu, \sigma, \delta, \delta')| \\
&\leq C(\underline{\mu}, \overline{\mu}, \underline{\sigma}, \overline{\sigma}, \nu) 2^{2J} (\Delta_\mu \vee \Delta_\sigma).
\end{aligned}
$$

(c) Control of $\sup_{(\delta, \delta') \in \{0, \ldots, N-1\}^2} \sup_{(\mu, \sigma) \in A_{\delta, \delta'}} |\Gamma_m^{(3)}(\mu, \sigma, \delta, \delta')|$. Assuming that $(h_1)$ holds, we have that for any $(\mu, \sigma)$ in $A_{\delta, \delta'}$, $x$ in the support of $f$,

$$
(5.24) \quad
\begin{aligned}
\left| f_0 \left( \frac{x - \mu}{\sigma} \right) \right. &- \left. f_0 \left( \frac{x - \mu_\delta}{\sigma_{\delta'}} \right) \right| \\
&\leq C_{f_0} \left[ |x| \left( \frac{\Delta_\sigma}{\underline{\sigma}^2} \right) + \frac{1}{\underline{\sigma}^2} (\overline{\sigma} \Delta_\mu + (|\overline{\mu}| \vee |\underline{\mu}|) \Delta_\sigma) \right].
\end{aligned}
$$

Hence, we derive from $(h_2)$ that there exists some positive constant $C(\underline{\mu}, \overline{\mu}, \underline{\sigma}, \overline{\sigma}, C_{f_0}, \nu)$ such that

$$
\sup_{(\delta,\delta')\in\{0,\ldots,N-1\}^2} \sup_{(\mu,\sigma)\in A_{\delta,\delta'}} |\Gamma_m^{(3)}(\mu,\sigma,\delta,\delta')|
$$
(5.25)
$$
\leq C(\underline{\mu}, \overline{\mu}, \underline{\sigma}, \overline{\sigma}, C_{f_0}, \nu)(\Delta_\mu \vee \Delta_\sigma).
$$

(d) Control of $\sup_{(\delta,\delta')\in\{0,\ldots,N-1\}^2} \sup_{(\mu,\sigma)\in A_{\delta,\delta'}} |\Gamma_m^{(4)}(\mu,\sigma,\delta,\delta')|$. It follows from (5.24) that there exists $C(\underline{\mu}, \overline{\mu}, \underline{\sigma}, \overline{\sigma}, C_{f_0}) > 0$ such that

$$
\sup_{(\mu,\sigma)\in A_{\delta,\delta'}} |\Gamma_m^{(4)}(\mu,\sigma,\delta,\delta')|
$$
$$
\leq C(\underline{\mu}, \overline{\mu}, \underline{\sigma}, \overline{\sigma}, C_{f_0})\left[\frac{1}{n}\sum_{i=1}^n (|X_i| + \underline{\mathbb{E}}(|X_i|) + 1)(\Delta_\mu \vee \Delta_\sigma)\right].
$$

Using again (5.21), we prove that there exists some positive constant $C(\underline{\mu}, \overline{\mu}, \underline{\sigma}, \overline{\sigma}, \nu, c, \beta, C_{f_0})$ such that

$$
\mathbb{P}_f\Bigg( \sup_{(\delta,\delta')\in\{0,\ldots,N-1\}^2} \sup_{(\mu,\sigma)\in A_{\delta,\delta'}} |\Gamma_m^{(4)}(\mu,\sigma,\delta,\delta')|
$$
(5.26)
$$
> C(\underline{\mu}, \overline{\mu}, \underline{\sigma}, \overline{\sigma}, \nu, c, \beta, C_{f_0})(\Delta_\mu \vee \Delta_\sigma)\log N \Bigg) \leq \frac{\beta}{4}.
$$

Collecting (5.22), (5.23), (5.25) and (5.26), we get

$$
\mathbb{P}_f\Bigg( \sup_{(\delta,\delta')\in\{0,\ldots,N-1\}^2} \sup_{(\mu,\sigma)\in A_{\delta,\delta'}} (\Gamma_m(\mu,\sigma) - \Gamma_m(\mu_\delta,\sigma_{\delta'}))
$$
(5.27)
$$
> C2^{2J}(\Delta_\mu \vee \Delta_\sigma)\log N \Bigg) \leq \frac{\beta}{2},
$$

for some constant $C = C(\underline{\mu}, \overline{\mu}, \underline{\sigma}, \overline{\sigma}, C_{f_0}, \|f\|_\infty, \nu, c, \beta, \varepsilon) > 0$. Finally, by setting $\Delta_\mu = \Delta_\sigma = n^{-2}2^{-J}$, we have that

$$
n^2 2^J(\overline{\mu} - \underline{\mu}) \wedge (\overline{\sigma} - \underline{\sigma}) \leq N \leq n^2 2^J(\overline{\mu} - \underline{\mu}) \vee (\overline{\sigma} - \underline{\sigma}) + 1.
$$

Hence, we deduce from (5.19) and (5.27) that if $n^2(\overline{\mu} - \underline{\mu}) \wedge (\overline{\sigma} - \underline{\sigma}) \geq 2$ and $n \geq 2$, there exists a constant $C = C(\underline{\mu}, \overline{\mu}, \underline{\sigma}, \overline{\sigma}, \|f_0\|_\infty, C_{f_0}, \|f\|_\infty, \nu, c, \beta, \varepsilon)$ such that

$$
\tau_m(\beta) \leq C\left( \frac{2^{J/2}}{n}\sqrt{\log(n^2 2^J)} + \frac{2^J}{n^2}\log^2(n^2 2^J) + \frac{\log(n^2 2^J)}{n} \right).
$$

This concludes the proof of Theorem 2.

5.5. *Proof of Corollary* 3. Since $\tilde{q}_{(2,D),\alpha} = t_{(2,D)}(u_\alpha)$ with the same notation as in Section 2, we have that

$$\tilde{q}_{(2,D),\alpha} \leq \tilde{W}_D(\alpha),$$

where

$$\tilde{W}_D(\alpha) = \frac{C(\alpha)}{n} \bigg( \left( \|f_0\|_\infty + \sqrt{\|f_0\|_\infty} \, \right) \sqrt{D \log \log n}$$
$$+ \frac{D(\log \log n)^2}{n} + \|f_0\|_\infty (\log \log n) \bigg),$$

$C(\alpha)$ being a positive constant. This upper bound is obtained by replacing in the proof of Proposition 1 $\alpha_n$ by $\alpha/2(1 + \log_2 n^2)$ and $C'\|f_0\|_\infty \log(D+1)$ by $C'\|f_0\|_\infty$ in (5.11).

Let $f$ belong to $\tilde{B}_s(R, M)$ and satisfy $(h_1)$ and $(h_2)$.

From Theorem 2 and the above upper bound, we deduce that there exists some positive constant $C = C(\underline{\mu}, \overline{\mu}, \underline{\sigma}, \overline{\sigma}, C_{f_0}, \|f_0\|_\infty, M, \nu, c, \alpha, \beta, s)$ such that if $f$ satisfies

$$\inf_{(\mu,\sigma)\in K} \left\| f - \frac{1}{\sigma} f_0\left( \frac{\cdot - \mu}{\sigma} \right) \right\|_2^2 \geq C \inf_{D\in\mathcal{D}_2} \bigg\{ R^2 D^{-2s} + \frac{\sqrt{D \log(n^2 D)}}{n}$$
$$+ \frac{D \log^2(n^2 D)}{n^2} + \frac{\log(n^2 D)}{n} \bigg\},$$

then

$$\mathbb{P}_f(\tilde{T}_\alpha \leq 0) \leq \beta.$$

Since $\mathcal{D}_2 = \{2^J, 0 \leq J \leq 2^{[\log_2(n^2/\log^3 n)]}\}$ and $n \geq 3$, there exists $c > 0$ such that

$$\inf_{D\in\mathcal{D}_2} \bigg\{ R^2 D^{-2s} + \frac{\sqrt{D \log(n^2 D)}}{n} + \frac{D \log^2(n^2 D)}{n^2} + \frac{\log(n^2 D)}{n} \bigg\}$$
$$\leq c \bigg( \inf_{D\in\mathcal{D}_2} \bigg\{ R^2 D^{-2s} + \frac{\sqrt{D \log n}}{n} \bigg\} + \frac{\log n}{n} \bigg),$$

and Corollary 3 can be proved in the same way as Corollary 1.

**Acknowledgment.** The authors want to express their thanks to the referees for their constructive suggestions.

LABORATOIRE DE STATISTIQUE
UNIVERSITÉ RENNES II
PLACE DU RECTEUR H. LE MOAL
CS 24307
35043 RENNES CEDEX
FRANCE
E-MAIL: magalie.fromont@uhb.fr

DÉPARTEMENT DE GÉNIE MATHÉMATIQUE
INSA
135 AVENUE DE RANGUEIL
31077 TOULOUSE CEDEX 4
FRANCE
E-MAIL: beatrice.laurent@insa-toulouse.fr